# Choiceless Polynomial Time

Andreas Blass, Yuri Gurevich and Saharon Shelah

**Abstract**

Turing machines define polynomial time (PTime) on strings but cannot deal with structures like graphs directly, and there is no known, easily computable string encoding of isomorphism classes of structures. Is there a computation model whose machines do not distinguish between isomorphic structures and compute exactly PTime properties? This question can be recast as follows: Does there exists a logic that captures polynomial time (without presuming the presence of a linear order)? Earlier, one of us conjectured the negative answer; see [74]. The problem motivated a quest for stronger and stronger PTime logics. All these logics avoid arbitrary choice. Here we attempt to capture the choiceless fragment of PTime. Our computation model is a version of abstract state machines (formerly called evolving algebras). The idea is to replace arbitrary choice with parallel execution. The resulting logic is more expressive than other PTime logics in the literature. A more difficult theorem shows that the logic does not capture all PTime.

## 1 Introduction

The standard computation model is Turing machines, whose inputs are strings. However, in combinatorics, database theory, etc., inputs are naturally structures (graphs, databases, etc.) indistinguishable up to isomorphism. In such cases, there is a problem with string presentation of input objects: there is no known, easily computable string encoding of isomorphism classes of structures. This calls for a computation model that deals with structures directly rather than via string encoding. There is a several such computation models in the literature, in particular relational machines [Abiteboul and Vianu 1991] and abstract state machines (formerly called evolving algebras) [Gurevich 1995].

The natural question is whether there is a computation model that captures PTime over structures (rather than strings). In different terms, essentially the same question has been raised in [Chandra and Harel 1982]. Gurevich translated Chandra–Harel's question as a question of existence of a logic that captures PTime and conjectured that no such logic exists [Gurevich 1988]. We address this issue in Section 3; here it suffices to say that the notion of logic is a very broad one and includes computation models.

If one seriously entertains the possibility that there is no logic that captures PTime, the question arises how much of PTime can be captured by a logic. Here we define a natural fragment of PTime and capture it by means of a version of abstract state machines (ASMs). We call the fragment Choiceless Polynomial Time (ČPTime). The idea is to eliminate arbitrary choice by means of parallel execution. Consider for example the Graph Reachability problem:

**Instance:** A graph $G = (V, E)$ with distinguished nodes $s$ and $t$ (an allusion to Source and Target respectively).

**Question:** Does there exists a path from $s$ to $t$ in $G$?

A common reachability algorithm constructs the set $X$ of all vertices reachable from $s$ and then checks if $X$ contains $t$. To construct $X$, an auxiliary "border-set" $Y \subseteq X$ is used.

```
if Mode = Initial then
    X, Y := {s}, Mode := Construct
endif

if Mode = Construct then
    if Y ≠ ∅ then
        choose y ∈ Y
            let Z = {z ∈ V − X : yEz}
                do in parallel
                    X := X ∪ Z
                    Y := (Y − {y}) ∪ Z
                enddo
            endlet
        endchoose
    else Mode := Examine
    endif
endif

if Mode = Examine then
    If t ∈ X then Output := Yes else Output := No endif
    Mode := Final
endif
```

If a given graph $G$ comes with an order on vertices, the order can be used to eliminate choice, but we are interested in structures which are not necessarily ordered. In the case of the reachability problem, choice can be eliminated by means of parallelism. Here is a revised version of the second transition rule.

```
if Mode = Construct then
   if Y ≠ ∅ then
      let Z = {z ∈ V − X : (∃y ∈ Y) yEz}
          X := X ∪ Z
          Y := Z
   else Mode := Examine
   endif
endif
```

Of course, one is not always so lucky. In Section 10, we describe a well-known PTime algorithm for the Perfect Matching problem. The algorithm uses choice and, as far as we know, there is no choiceless PTime algorithm for Perfect Matching.

Our computation model is explained in Section 4 and our formalization $\mathcal{L}$ of Choiceless PTime is given in Section 5. $\mathcal{L}$ is a computation model, but it can be viewed as a (very generalized) logic. In Section 7, we show that this logic is more expressive than Abiteboul–Vianu's logic (whose "formulas" are relational machines). In Section 10, we show that $\mathcal{L}$ does not express the parity of a naked set or the perfect matchability of a bipartite graph.

One shortcoming of $\mathcal{L}$ is that it is naturally three-valued: some input structures are accepted, some are rejected, and some may be neither accepted nor rejected by a given machine. In Section 11, we augment our computation model with explicit knowledge of the cardinality of (the base set of) the input structure. The augmented logic $\mathcal{L}'$ is two-valued. Parity is easily expressible in $\mathcal{L}'$ but perfect matchability is not expressible.

The logic $\mathcal{L}$ can be easily extended with a counting function (see Section 4) and maybe it should. This would be a natural way to continue this investigation. Let us notice though that, unless there is a logic that captures PTime, there is no end to possible extensions of $\mathcal{L}$. Any PTime decidable problem can be converted to a quantifier and added to $\mathcal{L}$.

Finally, there are some open problems that we leave. We conjecture that $\tilde{C}$PTime enjoys some sort of zero-one or convergence law. We conjecture that the extension $\tilde{C}$PTime + Count of $\tilde{C}$PTime with the counting function is a proper fragment of PTime. Further, it is not difficult to concoct an artificial complete problem for $\tilde{C}$PTime, but it would be interesting to see a natural one.

# 2 Preliminaries

We recall various definitions and establish some terminology and notation. In this paper, vocabularies are finite.

## 2.1 Global Relations

We start with a convenient notion of global relation [Gurevich 1988].

A $k$-ary *global relation* of vocabulary $\Upsilon$ is a function $\rho$ such that

- the domain $\text{Dom}(\rho)$ consists of $\Upsilon$-structures and is closed under isomorphism,

- $\rho$ associates a $k$-ary relation $\rho_A$ with every structure $A \in \text{Dom}(\rho)$, and

- $\rho$ is *abstract* in the following sense: every isomorphism from a structure $A \in \text{Dom}(\rho)$ onto a structure $B$ is also an isomorphism from the structure $(A, \rho_A)$ onto the structure $(B, \rho_B)$.

Typically the domain of a global relation of vocabulary $\Upsilon$ is the class of all $\Upsilon$-structures or the class of all finite $\Upsilon$-structures. For example, every first-order formula $\varphi(v_1, \ldots, v_k)$ of vocabulary $\Upsilon$ with free variables as shown, denotes a $k$-ary global relation $\rho(v_1, \ldots, v_k)$ on all $\Upsilon$-structures.

## 2.2 The Least Fixed Point Logic FO+LFP

The least fixed point logic has been around for a long time [Moschovakis 1974]. It is especially popular in finite model theory [Ebbinghaus and Flum 1995]. The latter book contains all the facts that we need. For the reader's convenience (and to establish notation), we recall a few things.

**Syntax**  FO+LFP is obtained from first-order logic by means of the following additional formula-formation rule:

- Suppose that $\varphi(P, \bar{v})$ is a formula with a $k$-ary predicate variable $P$ and a $k$-tuple $\bar{v}$ of free individual variables. Further suppose that $P$ occurs only positively in $\varphi$. If $\bar{t}$ be a $k$-tuple of terms, then

$$\left[LFP_{P,\bar{v}}(\varphi(P,\bar{v}))\right](\bar{t})$$

  is a formula.

The vocabulary and the free variables of the new formula $\psi$ are defined in the obvious way. In particular, $P$ is not in the vocabulary of $\psi$. Suppose that a predicate $Q$ belongs to the vocabulary of $\psi$; then $Q$ occurs only positively in $\psi$ if and only if it occurs only positively in $\varphi$.

**Semantics**  The formula $\varphi(P, \bar{v})$, may have free individual variables $\bar{u}$ in addition to $\bar{v}$. Let $\bar{w}$ be a $k$-tuple of fresh individual variables, and let $\Upsilon$ be the vocabulary of the formula

$$\psi(\bar{u}, \bar{w}) = \big[ LFP_{P,\bar{v}}(\varphi(P, \bar{u}, \bar{v})) \big](\bar{w})$$

The meaning of $\psi$ is a global relation $\rho(\bar{u}, \bar{w})$ whose domain consists of all $\Upsilon$-structures (unless we restrict the domain explicitly, for example to finite $\Upsilon$-structures).

Given an $\Upsilon$-structure $A$ with fixed values $\bar{a}$ of parameters $\bar{u}$, consider the following operator on $k$-ary relations over BaseSet($A$):

$$\theta(P) = \{\bar{v} : \varphi(P, \bar{a}, \bar{v})\}.$$

Since $\varphi$ is positive in $P$, $\theta$ is monotone in $P$. The $k$-ary relation $\rho_A(\bar{a}, \bar{w})$ is the least fixed point of $\theta$. To obtain the least fixed point, generate the sequence

$$\emptyset \subseteq \theta(\emptyset) \subseteq \theta^2(\emptyset) \subseteq \ldots$$

of $k$-ary relations over $A$. For a finite structure $A$, there exists a natural number $l$ such that $\theta^l(\emptyset) = \theta^{l+1}(\emptyset)$; in this case the least fixed point is $\theta^l(\emptyset)$. The case of infinite $A$ is similar except that $l$ may be an infinite ordinal.

**Simultaneous Induction**  Let $\Upsilon$ be a vocabulary and consider FO+LFP formulas $\varphi(P, Q, \bar{u})$ and $\psi(P, Q, \bar{v})$ of vocabulary $\Upsilon \cup \{P, Q\}$ which are positive in $P$ and $Q$. Here Arity($P$) = Length($\bar{u}$) and Arity($Q$) = Length($\bar{v}$). There may be additional free individual variables which we consider as paramenters.

Given an $\Upsilon$-structure $A$ with fixed parameters, consider the monotone operator

$$\theta(P, Q) = \Big( P \cup \{\bar{u} : \varphi(P, Q, \bar{u})\}, Q \cup \{\bar{v} : \psi(P, Q, \bar{v})\} \Big)$$

and let $(P^*, Q^*)$ be the least fixed point of $\theta$.

**Proposition 1**  *The global relations $P^*$ and $Q^*$ are expressible by FO+LFP formulas. Moreover, the result generalizes to simultaneous induction over any finite number of predicate variables.*

## 2.3 The Finite Variable Infinitary Logic

Again, the book [Ebbinghaus and Flum 1995] contains all the information that we need, but we recall a few things for the reader's convenience.

**Syntax** As in a popular version of first-order logic, $L^\omega_{\infty,\omega}$ formulas are built from atomic formulas by means of negations, conjunctions, disjunctions, the existential quantifier and the universal quantifier. The only difference is that $L^\omega_{\infty,\omega}$ allows one to form the conjunction and the disjunction of an arbitrary set $S$ of formulas provided that the total number of variables in all $S$-formulas is finite. Here we have an additional constraint: the vocabulary of every formula is finite. For every natural number $m$, $L^m_{\infty,\omega}$ is the fragment of $L^\omega_{\infty,\omega}$ where formulas use at most $m$ individual variables.

**Semantics** Every $k$-ary $L^\omega_{\infty,\omega}$ formula of vocabulary $\Upsilon$ denotes a $k$-ary global relation of vocabulary $\Upsilon$ in the obvious way.

An important fact is that every global relation expressible in FO+LFP is expressible in $L^\omega_{\infty,\omega}$ [Ebbinghaus and Flum 1995, Theorem 7.4.2].

**Pebble Games** There is a pebble game $G_m(A, B)$ appropriate to $L^m_{\infty,\omega}$. Here $A$ and $B$ are structures of the same purely relational vocabulary. For explanatory purposes, we pretend that that $A$ is located on the left and $B$ is located on the right, but in fact $A$ and $B$ may be the same structure.

The game is played by Spoiler and Duplicator. For each $i = 1, \ldots, k$, there are two pebbles marked by $i$: the left $i$-pebble and the right $i$-pebble. Initially all the pebbles are off the board. After any number of rounds, for every $i$, either both $i$-pebbles are off the board or else the left $i$-pebble covers an element of $A$ and the right $i$-pebble covers an element of $B$. In the obvious way, the pebbles on the board define a relation $R$ between $A$ to $B$. A round of $G^k(A, B)$ is played as follows.

If $R$ is not a partial isomorphism, then the game is over; Spoiler has won and Duplicator has lost. Otherwise Spoiler chooses a number $i$; if the $i$-pebbles are on the board, they are taken off the board. Then Spoiler chooses *left* or *right* and puts that $i$-pebble on an element of the corresponding structure. Duplicator puts the other $i$-pebble on an element of the other structure.

Duplicator wins a play of the game if the number of rounds in the play is infinite.

**Proposition 2** *If Duplicator has a winning strategy in $G_m(A, B)$, then no $L^m_{\infty,\omega}$ sentence distinguishes between $A$ and $B$. Therefore, for every FO+LFP sentence $\varphi$, there exists $m$ such that, for any $A$ and $B$, if Duplicator has a winning strategy in $G_m(A, B)$ then $\varphi$ does not distinguish $A$ from $B$.*

## 2.4 Set Theory

Let $A$ be a structure. In the literature, notation $|A|$ is used in two ways: to denote the base set of $A$ and to denote the cardinality of BaseSet($A$). We will employ notation $|A|$ only in the sense of cardinality; we will also use an alternative notation Card($A$) for the cardinality of $A$.

As usual in set theory, we identify a natural number (that is a non-negative integer) $i$ with the set of smaller natural numbers $\{j : j < i\}$; this set is called the von Neumann ordinal for $i$. The first infinite ordinal is denoted $\omega$.

We consider sets built from atoms (or urelements). The term *object* will mean an atom or a set. A set $X$ is *transitive* if $y \in x \in X$ implies $y \in X$. If $X$ is an object, then TC($X$) is the least transitive set $Y$ with $X \in Y$. An object $X$ is *hereditarily finite* if TC($X$) is finite.

$\mathbf{P}$ is the powerset operation; if $X$ is set then $\mathbf{P}(X)$ is the collection of all subsets of $X$. If $X$ is a finite set of atoms, then

$$\text{HF}(X) := \bigcup \{\mathbf{P}^i(X) : i < \omega\} = X \cup \mathbf{P}(X) \cup \mathbf{P}(\mathbf{P}(X)) \cup \ldots$$

Alternatively, HF($X$) can be defined as the smallest set $Y$ such that $X \subseteq Y$ and every finite subset of $Y$ is a member of $Y$.

Every set has an ordinal rank. If $x$ is an atom or the empty set, then the rank of $x$ equals 0; and if $x$ is a nonempty set, then the rank of $x$ equals 1 plus the maximum of the ranks of the elements of $x$.

# 3 PTime and PTime Logics

In this section, structures are finite and global relations are restricted to finite structures.

By definition, the PTime complexity class consists of languages, that is sets of (without loss of generality, binary) strings. A language $X$ is PTime if there exists a PTime Turing machine (that is polynomial time bounded Turing machine) that accepts exactly the strings in $X$. This definition is easily generalized to ordered structures by means of a standard encoding; see for example [Ebbinghaus and Flum 1995] or [Gurevich 1988]. We will say that a Turing machine accepts an ordered structure $A$ if it accepts the standard encoding of $A$.

The generalization to arbitrary (that is not necessarily ordered) structures is less obvious. One does not want to distinguish between isomorphic structures and there is no obvious, easily computable string encoding of isomorphism classes of structures.

The problem was first addressed by Chandra and Harel in the context of database theory [Chandra and Harel 1982]. We describe their approach. A *database* is defined as a purely relational structure whose elements come from some fixed countable set, without loss of generality the set of natural numbers. A *query* is a global relation over databases; recall that global relations respect isomorphisms. A query $Q$ is PTime if the set

$$\{(B, \bar{x}) : \bar{x} \in Q(B)\}$$

is PTime. Thus each PTime query $Q$ is given by a PTime Turing machine $M$ that accepts a string $s$ if and only if $s$ is the standard encoding of some $(B, \bar{x})$ with $\bar{x} \in Q(B)$; call $M$ a *PTime witness* for $Q$. Let $W$ be the collection of all PTime witnesses for all queries. It is easy to check that $W$ is not recursive. Chandra and Harel posed the following question. Does there exists a recursive set $S \subset W$ such that every PTime query has a PTime witness in $S$?

Gurevich translated their question as a question of the existence of a logic that captures PTime [Gurevich 1988]. He conjectured that the answer is negative; in this connection his definition of a logic is very broad. If desired, some obvious requirements can be imposed; see [Ebbinghaus 1985] in this connection. Here we recall Gurevich's definitions and slightly generalize them in order to define three-valued logics.

**PTime Global Relations**  What does it mean that a global relation is PTime? The question easily reduces to the case of nullary global relations. Indeed, let $\rho$ be a $k$-ary global relation of some vocabulary $\Upsilon$ and let $c_1, \ldots, c_k$ be the first $k$ individual constants outside of $\Upsilon$. Define the nullary relation $\sigma$ of vocabulary $\Upsilon^+ = \Upsilon \cup \{c_1, \ldots, c_k\}$ as follows. If $A$ is an $\Upsilon$-structure, and $a_1, \ldots, a_k$ are elements of $A$, and $B$ is the $\Upsilon^+$-expansion of $A$ where $a_1, \ldots, a_k$ interpret $c_1, \ldots, c_k$, then $\sigma_B \iff \rho_A(a_1, \ldots, a_k)$. Declare $\rho$ PTime if $\sigma$ is so.

A nullary global relation $\rho$ of vocabulary $\Upsilon$ can be identified with the class of $\Upsilon$-structures $A$ such that $\rho_A$ is true. It remains to define what does it mean that a class $K$ of structures of some vocabulary $\Upsilon$ is PTime. Let $<$ be a binary predicate not in $\Upsilon$. An *ordered version* of an $\Upsilon$-structure $A$ is a structure $B$ of vocabulary $\Upsilon \cup \{<\}$ such that the $\Upsilon$-reduct of $B$ is isomorphic to $A$ and the interpretation of $<$ is a linear order.

Define a class of $K$ of $\Upsilon$-structures to be PTime if it is closed under isomorphisms and there exists a PTime Turing machine $M$ (a PTime witness for $K$) which accepts a binary string $s$ if and only if $s$ is the standard encoding of an ordered version of some structure in $K$.

**Logics**  For simplicity, we define logics whose formulas denote nullary global relations. The trick above allows one to extend such a logic so that its formulas denote arbitrary global relations.

A *logic* $L$ is given by a pair of functions (Sen,Sat) satisfying the following conditions. Sen associates with every vocabulary $\Upsilon$ a recursive set $\text{Sen}(\Upsilon)$ whose elements are called *L-sentences of vocabulary* $\Upsilon$. Sat associates with every vocabulary $\Upsilon$ a recursive relation $\text{Sat}_\Upsilon(A, \varphi)$ where $A$ is an $\Upsilon$-structure and $\varphi$ an $L$-sentence of vocabulary $\Upsilon$. We say that $A$ satisfies $\varphi$ (symbolically $A \models \varphi$) if $\text{Sat}_\Upsilon(A, \varphi)$ holds. It is assumed that $\text{Sat}_\Upsilon(A, \varphi) \iff \text{Sat}_\Upsilon(B, \varphi)$ if $A$ and $B$ are isomorphic.

If $\varphi$ is a sentence of vocabulary $\Upsilon$, let $\text{Mod}(\varphi)$ be the collection of $\Upsilon$-structures $A$ satisfying $\varphi$.

**PTime Logics**  Let $L$ be a logic. For each $\Upsilon$ and each $\varphi \in \text{Sen}(\Upsilon)$, let $K(\Upsilon, \varphi)$ be the class of $\Upsilon$-structures $A$ such that $A \models \varphi$. Call $L$ PTime, if every class $K(\Upsilon, \varphi)$ is PTime.

**Logic that Capture Ptime**  A logic $L$ *captures PTime* if it is PTime and, for every vocabulary $\Upsilon$, every PTime class of $\Upsilon$-structures coincides with some $K(\Upsilon, \varphi)$.

**Remark**  It may seem odd that the definition of logic does not require any uniformity in $\Upsilon$, but uniformity is not necessary. It is not hard to show that, if there is a logic of graphs (rather than arbitrary structures) that captures PTime on graphs, then there is a logic that captures PTime on arbitrary structures and that possesses some uniformity with respect to $\Upsilon$. [Gurevich 1988].  □

**Three-Valued Logics**  In the cases when a logic is really a computation model and sentences are computing machines, $A$ satisfies $\varphi$ means that $\varphi$ accepts $A$. That calls for the following natural generalization. Call logics defined above *two-valued*.

A *three-valued logic* $L$ is like a two-valued logic except that $\text{Sat}_\Upsilon(A, \varphi)$ has three possible values telling us whether $\varphi$ accepts $A$, or $\varphi$ rejects $A$, or neither. It is assumed that $\text{Sat}_\Upsilon(A, \varphi) = \text{Sat}_\Upsilon(B, \varphi)$ if $A$ and $B$ are isomorphic.

Each $L$-sentence $\varphi$ of vocabulary $\Upsilon$ gives rise to two disjoint classes of $\Upsilon$-structures. The class $\text{Mod}^+(\varphi)$ of $\Upsilon$-structures accepted by $\varphi$ and the class $\text{Mod}^-(\varphi)$ of classes rejected by $\varphi$. Call $L$ PTime if, for every $\varphi$, the classes $\text{Mod}^+(\varphi)$ and $\text{Mod}^-(\varphi)$ are PTime.

Two disjoint classes $K_1, K_2$ of structures of some vocabulary $\Upsilon$ are *L-separable* if there exists an $L$-sentence $\varphi$ such that $K_1 \subseteq \text{Mod}^+(\varphi)$ and $K_2 \subseteq \text{Mod}^-(\varphi)$. We will see that this is a more robust notion than the similar notion where $\subseteq$ is replaced with equality.

**Abiteboul-Vianu Relational Machines**  Finally, for future reference, we recall (a version of) Abiteboul-Vianu's Relational machines [Abiteboul and Vianu 1991; Abiteboul-Vardi-Vianu 1997].

A relational machine is a Turing machine augmented with a *relational store* which is a structure of a fixed purely relational vocabulary $\Upsilon$. A part $\Upsilon_0$ of the vocabulary is devoted to input relations. The Turing tape is initially empty. As usual, the program consists of "if condition then action" instructions. Here is an example of an instruction.

> If the control state is $s_3$, and the head reads symbol 1, and the relation $R_1$ is empty, then change the state to $s_4$, replace 1 by 0, move the head to the right and replace $R_2$ with $R_2 \cap R_3$.

In general, instructions are Turing instructions except that (1) the condition may be augmented with the emptyness test of one of the relations, and (2) the action may be augmented with an algebraic operation on the relations. The algebraic operations are of the following four types. It is assumed that the arities of the operations involved are appropriate; in the example above, the relations $R_2$ and $R_3$ are all of the same arity.

- Boolean operations.

- Projections $\pi_{i_1...i_m} R_k$. Project $R_k$ on the coordinates $i_1, \ldots, i_m$ in the specified order.

- Cartesian product of two relations.

- Selections $\sigma_{i=j} R_k$. Select the tuples in $R_k$ whose $i$-th component coincides with the $j$-th component.

A PTime relational machine $M$ can be defined as a relational machine together with a polynomial $p(n)$ bounding the number of computation steps on input structure of size $n$. The notion of PTime relational machine gives rise to a PTime logic (whcih may be called AV Logic) where the sentences of vocabulary $\Upsilon_0$ are PTime relational machines with input vocabulary $\Upsilon_0$.

In our view, AV Logic is naturally three-valued. Given an input structure $I$ of size $n$, a PTime relational machine $(M, p(n))$ may accept $I$ within time $p(n)$, may reject $I$ within time $p(n)$, or do neither. One cannot modify $M$ so that it computes $p(n)$ because $n$ is unknown to the machine. Of course one can clump together all non-accepted structures but this definition is too sensitive to the choice of $p(n)$.

# 4 Computation Model

Our computing devices are abstract state machines (ASMs, formerly called evolving algebras) [Gurevich 1995, Gurevich 1997] adapted for our purposes here.

## 4.1 Vocabularies

An ASM vocabulary is a finite collection of function names, each of a fixed arity. Some function names may be marked as *relational* or *static*, or both. Relational names are also called *predicates*. A function name is *dynamic* if it is not marked static. The Greek letter $\Upsilon$ is reserved to denote vocabularies.

In our case, every vocabulary consists of the following four parts:

**Logic names** The equality sign, nullary function names *true*, *false* and the names of the usual Boolean operations. All logic names are relational and static. (The standard ASM definition [Gurevich 1995] requires another logic name, *undef*, but we will not employ *undef* here, using $\emptyset$ instead as a default value.)

**Set-theoretic names** The static binary predicate $\in$ and the following static non-predicate function names.

- Nullary names $\emptyset$ and Atoms.
- Unary names $\bigcup$ and TheUnique.
- A binary name Pair.

**Input names** A finite collection of static names. For simplicity of exposition, we assume that all input names are relational.

**Dynamic names** A finite collection of dynamic function names including nullary predicates Halt and Output.

## 4.2 States

A state $A$ of vocabulary $\Upsilon$ is a structure $A$ of vocabulary $\Upsilon$ satisfying a number of conditions described in this subsection.

**Base Set** The base set of $A$ consists of two disjoint parts:

1. A finite set $X$ of atoms, that is elements that are not sets.
2. The collection of all hereditarily finite sets built from the atoms.

The atoms and the sets are *objects* of $A$. The objects form a transitive set $\mathrm{HF}(X)$ which can be defined as the closure of $X$ under the following operation: If $n$ is a natural number and $x_1, \ldots, x_n$ are in, then throw $\{x_1, \ldots, x_n\}$ in. We have also:

$$\mathrm{HF}(X) = \bigcup_{n<\omega} \mathrm{PowerSet}^n(X).$$

**Set-Theoretic Functions**  The interpretations of $\in$ and $\emptyset$ are obvious. Atoms is the set of atoms. If $a$ is an atom or a logic element, then $\bigcup a = \emptyset$. If $a_1, \ldots, a_j$ are atoms and $b_1, \ldots, b_k$ are sets then $\bigcup \{a_1, \ldots, a_j, b_1, \ldots, b_k\} = b_1 \cup \cdots \cup b_k$. If $a$ is a singleton set, then TheUnique($a$) is the unique element of $a$; otherwise TheUnique($a$) = $\emptyset$. (If $x$ is a set then $x$ = TheUnique$\{x\}$ = $\bigcup\{x\}$, so TheUnique is redundant in this situation, but it is needed if $x$ is an atom.) Pair($a, b$) = $\{a, b\}$.

**Logic Names**  *false* and *true* are interpreted as 0 and 1 respectively. Recall that 0 is $\emptyset$ and that 1 is Pair(0,0) = $\{\emptyset\}$. The Boolean connectives are interpreted in the obvious way over the Boolean values $0, 1$ and take the value 0 if at least one of the arguments is not Boolean.

**Predicates**  Predicates are interpreted as functions whose only possible values are the Boolean values $0, 1$. If $P(\bar{a})$ evaluates to 1 (respectively 0), we say that $P(\bar{a})$ *holds* or *is true* (respectively, *fails* or *is false*). The input predicates "live" over the atoms: if $P$ is an input predicate and $P(a_1, \ldots, a_j)$ holds, then every $a_i$ is an atom.

**Dynamic Functions**  Define the *extent* of a dynamic function $f$ of arity $j$ to be the set

$$\{(x_0, \ldots, x_j) : f(x_0, \ldots, x_{j-1}) = x_j \neq 0\}.$$

The only restriction on the interpretation of a dynamic function $f$ is that its extent is finite.

## 4.3 Input Structures

Consider an ASM vocabulary $\Upsilon$. An input structure appropriate for $\Upsilon$-programs is any finite structure $I$ of the input vocabulary. A little problem arises in the case when some elements of $I$ are sets. The actual input corresponding to $I$ is a structure isomorphic to $I$ whose base set (the universe) consists of atoms. An $\Upsilon$-state is *initial* if the extent of every dynamic function is empty. For any input structure $I$ appropriate for $\Upsilon$, there is a unique, up to isomorphism, initial $\Upsilon$-state $A$ where the atoms together with input relations form a structure isomorphic to $I$.

Remark. We will be not very careful in distinguishing between an input structure $I$ and its atomic version. Without loss of generality, one may assume that the input structure itself consists of atoms.

## 4.4 Terms

By induction, we define a syntactic category of terms and a subcategory of Boolean terms.

- A variable is a term.

- If $f$ is a function name of arity $j$ and $t_1, \ldots, t_j$ are terms, then $f(t_1, \ldots, t_j)$ is a term. If $f$ is a predicate then $f(t_1, \ldots, t_j)$ is Boolean.

- Suppose that $v$ is a variable, $t(v)$ is a term, $r$ is a term without free occurrences of $v$, and $g(v)$ is a Boolean term. Then

$$\{t(v) : v \in r : g(v)\}$$

  is a term.

In the usual way, the same induction is used to defined free variables of a given term. In particular, the free variables of $\{t(v) : v \in r : g(v)\}$ are those of $t(v), r$ and $g(v)$ except for $v$.

Semantics is obvious. In particular, the value of $\{t(v) : v \in r : g(v)\}$ at a given state $A$ is the set of values $\text{Val}_A(t(v))$ such that, in $A$, $v \in r$ and $g(v)$ holds.

## 4.5 Syntax of Rules

Transition rules are defined inductively.

**Skip**  Skip is a rule.

**Update Rules**  Suppose that $f$ is a dynamic function name of some arity $r$ and $t_0, \ldots, t_r$ are terms. If $f$ is relational, we require that $t_0$ is Boolean. Then

$f(t_1, \ldots, t_r) := t_0$

is a rule.

**Conditional Rules**  If $g$ is a Boolean term and $R_1, R_2$ are rules, then

`if` $g$ `then` $R_1$ `else` $R_2$ `endif`

is a rule.

**Do-forall Rules** If $v$ is a variable, $r$ is a term without $v$ free, and $R_0(v)$ is a rule, then

```
do forall v ∈ r
    R₀(v)
enddo
```

is a rule with *head variable* $v$, *guard* $r$ and *body* $R_0$. The definition of free and bound variables is obvious.

Abbreviate rule

```
do forall v ∈ {0,1}                              do-in-parallel
    if v = 0 then R₀                to               R₀, R₁
    else R₁ endif                                    enddo
enddo
```

Here and in the following, we use standard notation as a more readable substitute for the official syntax. In particular, 0 means $\emptyset$, $\{x,y\}$ means $\mathrm{Pair}(x,y)$, $\{x\}$ means $\mathrm{Pair}(x,x)$, and 1 means $\{0\}$.

## 4.6 Semantics of Rules

If $\zeta$ is a variable assignment over a state $A$, assigning values to finitely many variables, then the pair $B = (A, \zeta)$ is an *expanded state*, $A = \mathrm{State}(B)$, $\zeta = \mathrm{Assign}(B)$, and $\mathrm{Dom}(\zeta) = \mathrm{Var}(B)$. Further, let $v$ be a variable and $a$ an element of $A$. Then $B(v \mapsto a)$ is the expanded state obtained from $B$ by assigning or reassigning $a$ to $v$. In other words, $B(v \mapsto a) = (A, \zeta')$ where $\mathrm{Dom}(\zeta') = \mathrm{Dom}(\zeta) \cup \{v\}$, $\zeta'(v) = a$ and $\zeta'(u) = \zeta(u)$ for the remaining variables.

A *location* of an expanded state $A$ is a pair $\ell = (f, \bar{a})$ where $f$ is a dynamic function and $\bar{a}$ is a tuple of elements of $A$ such that the length of $\bar{a}$ equals the arity of $f$. If $b$ is also an element of $A$, then the pair $\alpha = (\ell, b)$ is an *update* of $A$. $((f, (a_1, \ldots, a_j)), b)$ is abbreviated to $(f, a_1, \ldots, a_j, b)$. To *fire* $\alpha$ at $A$, put $b$ into the location $f$, that is, redefine $A$ so that $f(\bar{a}) = b$. The other locations remain intact. The resulting state is the *sequel* of $A$ with respect to $\alpha$. Two updates *clash* if they have the same locations but different new contents.

An *action* over a state $A$ is a set of updates of $A$. An action is *consistent* if it contains no clashing updates. To perform an action $\beta$, do the following. If $\beta$ is consistent, then fire all updates $\alpha \in \beta$ simultaneously; otherwise do nothing. The result is the *sequel* of $A$ with respect to $\beta$. If $\beta$ is inconsistent then the sequel of $A$ is $A$ itself.

A rule $R$ and the expanded state $A$ are *appropriate* for each other if $\text{Voc}(A)$ contains all function symbols in $R$ and $\text{Var}(A)$ contains all free variables in $R$.

Now we are ready to explain the semantics of rules. The denotation $\text{Den}(R)$ of a rule $R$ is a function on expanded states $A$ appropriate for $R$. Each $\text{Den}(R)(A)$ (or $\text{Den}(R, A)$ for brevity) is an action. To fire $R$ at $A$, perform the action $\text{Den}(R, A)$ at $\text{State}(A)$. The *sequel* of $A$ with respect to $R$ is the sequel of $A$ with respect to $\text{Den}(R, A)$. $\text{Den}(R, A)$ is defined by induction on $R$.

**Skip**  $\text{Den}(\text{Skip}, A) = \emptyset$.

**Update Rules**  If $R$ is an update rule $f(\bar{s}) := t$ and $\ell$ is location $(f, \text{Val}_A(\bar{s}))$, then $\text{Den}(R, A) = \{(\ell, \text{Val}_A(t))\}$.

**Conditional Rules**  If $R$ is the rule  `if g then` $R_1$ `else` $R_2$ `endif`, then

$$\text{Den}(R, A) = \begin{cases} \text{Den}(R_1, A) & \text{if } g \text{ holds at } A; \\ \text{Den}(R_2, A) & \text{otherwise.} \end{cases}$$

**Do-forall Rules**  If $R$ is

```
do forall v ∈ r,
   R_0(v)
enddo
```

then

$$\text{Den}(R, A) = \bigcup \{\text{Den}(R_0(v), A(v \mapsto a)) : a \in \text{Val}_A(r)\}$$

## 4.7  Programs

A program is a rule without free variables. The vocabulary $\text{Voc}(\Pi)$ of a program $\Pi$ is the collection of function names that occur in $\Pi$. States of $\Pi$ are states of the vocabulary $\text{Voc}(\Pi)$.

**Runs**  A *run* of $\Pi$ is a (finite or infinite) sequence $\langle A_i : i < \kappa \rangle$ of states of $\Pi$ such that

- $A_0$ is an initial state,

- every $A_{i+1}$ is a sequel of $A_i$ with respect to $\Pi$, and
- Halt fails at every $A_i$ with $i + 1 < \kappa$.

Here $\kappa$ is a positive integer or the first infinite ordinal $\omega$. The *length of a finite run* $\langle A_i : i \leq l \rangle$ is $l$. The *length of an infinite run* is $\omega$.

Let $I$ be an input structure for $\Pi$. *The run of* $\Pi$ *on* $I$ is the run $\langle A_i : i < \kappa \rangle$ such that

- $A_0$ is the initial state generated by $I$, and
- either $\kappa$ is infinite, or else $\kappa$ is finite and Halt holds at the final state $A_{\kappa-1}$.

The *base set* and *objects* of a run $\langle A_i : i < \kappa \rangle$ are those of $A_0$.

## 4.8 The Counting Function

There are many ways to extend the computation model described above without introducing explicit choice in its full generality. One gentle extension of the computation model described above is considered in Section 11. A more powerful extension is achieved by introducing the counting function which, given a set $x$ of cardinality $k$, produces the von Neumann ordinal for $k$.

Remark. Since the computation model is expandable by adding static functions for counting or perfect matching, etc., one gets in fact a notion of relative computability.

# 5 Choiceless PTime

It is easy to check that every computable global relation on finite structures is computable by an appropriate ASM program. But we are interested in polynomial time computations.

## 5.1 The Definition of Choiceless PTime

**Critical, Active and Relevant Objects**   Let $A$ be a state and $x \in \mathrm{BaseSet}(A)$.

- Object $x$ is *critical* at $A$ if $x$ is an atom, or $x \in \{0, 1\}$, or $x$ is a value of a dynamic function, or $x$ is a component of a location where some dynamic function takes a value different from $\emptyset$.
- Object $x$ is *active* at $A$ if $x \in \mathrm{TC}(y)$ for some critical $y$.

Further, let $\rho$ be a run of a program $\Pi$. An object $x$ of $\rho$ is *active in* $\rho$ if it is so at some state of $\rho$.

**PTime Programs** One may define a notion of honest computation time. Don't just count computation steps; charge for example a unit of time for every function application, every update, etc. In other words, instead of macrosteps, count the microsteps of the computation [BG]. Here we use a simpler way to bound the computation time.

A *PTime (bounded) program* $\bar{\Pi}$ is a triple $\bar{\Pi} = (\Pi, p(n), q(n))$ where $\Pi$ is a program and $p(n), q(n)$ are integer polynomials. The *run* of $\bar{\Pi}$ on an input structure $I$ of size $n$ is the longest initial segment $\rho$ of the run of $\Pi$ on $I$ such that the length of $\rho$ is $\leq p(n)$ and the number of active objects in $\rho$ is $\leq q(n)$. A PTime program $\bar{\Pi}$ *accepts* (respectively *rejects*) an input structure $I$ if the run of $\bar{\Pi}$ on $I$ halts and Output equals *true* (respectively *false*) in the final state.

**Choiceless Polynomial Time** Notice the classes of accepted and rejected input structures are disjoint but not necessarily complementary and that increasing the polynomial bound may increase these classes. If the size of the input structure is known and if one uses a version of honest computation time, a program may keep track of the honest computation time and in this way insure that every computation is halting. We will consider this scenario in Section 11.

Here, we define Choiceless Polynomial Time (in brief $\tilde{C}$PTime) as a collection of pairs $(K_1, K_2)$ where $K_1, K_2$ are disjoint classes of finite structures of the same vocabulary. A pair $(K_1, K_2)$ is in $\tilde{C}$PTime (or $\tilde{C}$PTime separable) if there exists a PTime program that accepts all structures in $K_1$ and rejects all structures in $K_2$. The program may accept some structures not in $K_1$ or reject some structures not in $K_2$. Obviously, there is a three-valued logic that expresses exactly $\tilde{C}$PTime properties; use PTime programs as sentences.

A class $K$ of finite structures of the same vocabulary $\Upsilon$ is in $\tilde{C}$PTime, if the pair $(K, K')$ is in $\tilde{C}$PTime where $K'$ is the complement of $K$ in the class of finite structures of vocabulary $\Upsilon$.

Call two programs $\Pi$ and $\Sigma$ *PTime equivalent* if

- for every PTime version $\bar{\Pi}$ of $\Pi$, there exists a PTime version $\bar{\Sigma}$ of $\Sigma$ which accepts all input structures accepted by $\bar{\Pi}$ and rejects all input structures rejected by $\bar{\Pi}$, and

- for every PTime version $\bar{\Sigma}$ of $\Sigma$, there exists a PTime version $\bar{\Pi}$ of $\Pi$ which accepts all input structures accepted by $\bar{\Sigma}$ and rejects all input structures rejected by $\bar{\Sigma}$.

## 5.2 Upper Bounds for ČPTime

In our definition of a PTime program, the given polynomials bound not only time but also space. We show that ČPTime is not too broad.

**Theorem 1** *Consider a PTime program $\bar{\Pi} = (\Pi, p(n), q(n))$.*

1. *There is a PTime-bounded Turing machine that accepts exactly those strings that encode ordered versions of input structures accepted by $\bar{\Pi}$ and rejects exactly those strings that encode ordered versions of input structures rejected by $\bar{\Pi}$.*

2. *There exists a polynomial $r(n)$ such that the number of microsteps in every run of $\bar{\Pi}$ on an input structure of size $n$ is bounded by $r(n)$.*

**Proof** 1. The desired Turing machine simulates the given PTime program. The bound $r$ in a term $\{s(v) : v \in r : g(v)\}$ and in a do-forall rule ensures that the number of immediate subcomputations is bounded by the number of active elements and thus by $q(n)$. This yields a polynomial bound on the work needed to simulate one transition in the run. since the number of transitions is bounded by $p(n)$, the whole simulation takes only polynomial time.

2. Since the number of macrosteps is bounded by a polynomial, it suffices to check that the number of microsteps needed to fire an arbitrary rule $R$ is bounded by a polynomial. This is done by an obvious induction on $R$.  □

The part 1 of the theorem gives the following corollary.

**Corollary 1** *Every ČPTime pair of structure classes $(K_1, K_2)$ is separated by a PTime class.*

## 5.3 A Lower Bound for CP

In the previous subsection, we have shown that our definition of CP is not too broad. One may also worry that it is too narrow, that — because of the use of transitive closure in the definition of active objects — it is possible to create a large number of active objects in short time. Some active objects, namely atoms and $0, 1$, exist already in the initial state. The problem is to show that only polynomially many active nonempty sets can be created within polynomial honest computation time. We show that, under the definition of honest computation time hinted at above, the number of objects activated (that is the number of active objects which are inactive in the initial state) in any run of a PTime program is bounded by the honest computation time. The details of the definition of honest computation time are not important.

Consider a PTime program $\Pi$. Without loss of generality, we may assume that $\Pi$ does not reuse variables, that is no variable is bound more than once. It follows that, in every subrule of $\Pi$, no variable is bound more than once and no variable occurs both free and bound. Define a *grounded term* to be a pair $(t, A)$ where $t$ is a term and $A$ is an expanded state appropriate for $t$. Similarly, define a *grounded rule* to be a pair consisting of a rule and expanded state appropriate for it. Notice that a grounded term $(t, A)$ has a value, namely $\mathrm{Val}_A(t)$.

The following definitions are intended to describe, for each grounded term or rule, say $(X, A)$, a partially ordered set (poset) $\mathrm{Pre}(X, A)$ whose nodes are labeled with grounded terms that one would naturally evaluate in the course of evaluating $X$ at $A$; the order of $\mathrm{Pre}(X, A)$ reflects the order in which one would evaluate the grounded terms. $\mathrm{Pre}(X, A)$ is similar to the parse tree of $X$, but there are some distinctions. To prevent the definitions from getting even longer that they are, we omit the grounded terms involved in evaluating guards; one could include them without any damage to our argument. Also, we omit subrules that end up not being used because of clashes.

In fact, $\mathrm{Pre}(X, A)$ is not necessarily a tree. It will be convenient for our purposes that, for each free variable of $X$, there is at most one node with a label of the form $(x, A)$ or $(x, A(\bar{v} \mapsto \bar{a}))$; this gives rise to the following auxiliary definition. Let $P$ be a poset whose nodes are labeled with ground terms, and let $F$ be a collection of variables such that no label in $P$ has the form $(x, A(x \mapsto a))$ and each node with a label of the form $(x, A)$ or $(x, A(\bar{v} \mapsto \bar{a}))$ is minimal in $P$. Then *adjusting $P$ with respect to $F$* means merging, for each $x \in F$, all nodes of $P$ with labels of the form $(x, A)$ or $(x, A(v \mapsto a))$ into one node labeled with $(x, A)$.

Define a *disjoint union of labeled posets* in the obvious way: order and the labels within each piece are preserved and elements of distinct pieces are incomparable. Now we are ready define posets $\mathrm{Pre}(X, A)$ by induction on $X$.

**Definition 1** $\mathrm{Pre}(t, A)$ is defined by recursion on $t$.

- If $t$ is a variable $x$, then $\mathrm{Pre}(t, A)$ is a singleton poset whose only node is labeled with $(x, A)$.

- If $t$ is $f(t_1, \ldots, t_j)$, then $\mathrm{Pre}(t, A)$ is obtained from the disjoint union of $\mathrm{Pre}(t_1, A), \ldots, \mathrm{Pre}(t_j, A)$ by adding a $(t, A)$-labeled node at the top and adjusting the result with respect to the free variables of $t$.

- If $t$ is $\{s(v) : v \in r : g(v)\}$, then $\mathrm{Pre}(t, A)$ is obtained from the disjoint union of labeled posets $\mathrm{Pre}(s(v), A(v \mapsto a))$, for all $a \in \mathrm{Val}_A(r)$ such that $\mathrm{Val}_A(g(a)) = true$, by doing the following. (1) Adjoin a $(t, A)$-labeled node at the top. (2) Adjoin a copy of $\mathrm{Pre}(r, A)$ below all $(v, A(v \mapsto a))$-labeled nodes. (3) Adjust the result with respect to the free variables of $t$.

□

Note that Pre($t, A$) always has the top node labeled with ($t, A$). Further, for each free variable $x$ of $t$, there is at most one node labeled with ($x, A$) and this node (if present at all) is minimal in Pre($t, A$).

**Definition 2** Pre($R, A$) is defined by recursion on $R$.

- If $R$ is Skip, then Pre($R, A$) = $\emptyset$.

- If $R$ is $f(t_1, \ldots, t_j) := t_0$, then construct Pre($R, A$) as follows. Form the disjoint union of Pre($t_0, A$), ..., Pre($t_j, A$) and adjust the result with respect to the free variables of $R$.

- If $R$ is "if $g$ then $R_1$ else $R_2$ endif" then Pre($R, A$) is Pre($R_1, A$) or Pre($R_2, A$) according to whether Val$_A(g)$ is true or false.

- If $R$ is "do forall $v \in r$, $R_0(v)$ enddo", then construct Pre($R, A$) as follows. Form the disjoint union of Pre($R_0(v), A(v \mapsto a)$) for all $a \in$ Val$_A(r)$. Add a copy of Pre($r, A$) below each ($v, A(v \mapsto a)$) if there any; otherwise add a copy of Pre($r, A$) to the disjoint union. Adjust the result with respect to the set of free variables of $R$.

□

If ($X, A$) is a grounded term or rule, let Val[Pre($X, A$)] be the collection of objects Val$_B(s)$ such that ($s, B$) is a label in Pre($X$).

**Lemma 1**  1. If $Den(R, A)$ contains an update $(f, (a_0, \ldots, a_{j-1}), a_j)$, then every $a_i \in Val[Pre(R, A)]$.

2. Suppose that ($X, A$) is a grounded term or rule with bound variable $v$. If ($v, B$) is a label in $Pre(X, A)$, then $B$ has the form $C(v \mapsto a)$ where $C = A(\bar{u} \mapsto \bar{b})$ and the variables $\bar{u}$ (if present at all) are all different from $v$.

**Proof**
1. Induction on $R$.
2. Induction on $X$.  □

The labels of Pre($X, A$) are (some of the) grounded terms that would be evaluated when one evaluates $X$ in $A$. At least one unit of honest computation time should be spent to evaluate each of the labels. So if a run $\langle A_0, \ldots, A_l \rangle$ of the program $\Pi$ takes honest computation time $T$, then

$$\mathrm{Card}(\bigcup_i \mathrm{Val}[\mathrm{Pre}(\Pi, A_i)]) \leq T.$$

**Theorem 2** *Consider a run $\rho = \langle A_0, \ldots, A_l \rangle$ of $\Pi$. Every object $x$ activated in $\rho$ belongs to $\bigcup_i \text{Val}[\text{Pre}(\Pi, A_i)]$.*

**Proof** Call sets $0, 1$ binary and let $x$ be an active nonbinary set in $\rho$. In view of the convention about dynamic functions in initial states, $A_0$ has no critical nonbinary sets and thus no active nonbinary sets. Let $i$ be the first index such that $x$ is active in $A_{i+1}$. So there is a nonbinary set $y$, critical for $A_{i+1}$, with $x \in \text{TC}(y)$. Since $y$ is not critical for $A_i$, there must be an update, executed in the step from $A_i$ to $A_{i+1}$, involving $y$ as either the new value or a component of the location. By Lemma 1, $y \in \text{Val}[\text{Pre}(\Pi, A_i)]$. Among all nodes $n$ in $\text{Pre}(\Pi, A_i)$ such that $x \in \text{TC}(\text{Val}(\text{Label}(n)))$, choose a minimal one. Call this node $n_0$ and let $(t, B) = \text{Label}(n_0)$. Our goal is to show that $\text{Val}_B(t) = x$. So suppose this fails. Then there exists $w \in \text{Val}_B(t)$ such that $x \in \text{TC}(w)$. (This includes the possibility that $x = w$.) As $x$ is a nonbinary set, $w$ is a nonbinary set. We consider the various possibilities for $t$ and deduce a contradiction in every case. Note that $B$ is an expansion of $A_i$.

Suppose that $t$ has the form $f(\bar{s})$ for a dynamic $f$. Then $\text{Val}_B(t)$ is critical already in $A_i$ and therefore $x$ is active in $A_i$, contrary to our choice of $i$.

Suppose that $t$ is $\emptyset$ or Atoms. This is absurd, as $\text{Val}_B(t)$ contains a nonbinary set $w$.

Suppose that $t$ is $\bigcup s$. Since $w \in \text{Val}_B(t) = \bigcup \text{Val}_B(s)$, we have $w \in u \in \text{Val}_B(s)$ for some $u$. Since $x \in \text{TC}(w)$, we have $x \in \text{TC}(\text{Val}_B(s))$. But $\text{Pre}(t, B)$ has a node labeled $(s, B)$ and thus there is a node labeled $(s, B)$ below $n_0$ in $\text{Pre}(\Pi, A)$. This contradicts the choice of $n_0$.

Suppose $t$ is $\{s_1, s_2\}$. Then $x \in \text{TC}(w) = \text{TC}(\text{Val}_B(s))$ for some $s \in \{s_1, s_2\}$. The rest is as in the $\bigcup$ case.

Suppose that $t$ is $\text{TheUnique}(s)$. Since $\text{Val}_B(t)$ is a set (not an atom), $\text{TheUnique}(s) = \bigcup(s)$ here, and we get a contradiction as in the $\bigcup$ case.

Suppose that $t$ is $P(\bar{s})$ where $P$ is a predicate name. According to our presentation of truth values, $\text{Val}_B(t)$ is either $\emptyset$ or $\{\emptyset\}$. In the first case, we get a contradiction as in the $\emptyset$ case. In the second case, $w = \emptyset$ which is impossible as $\text{TC}(w)$ contains a nonbinary set $x$.

Suppose that $t$ is $\{s(v) : v \in r : g(v)\}$. Since $w \in \text{Val}_B(t)$, there is some $a \in \text{Val}_B(r)$ such that $\text{Val}_{B(v \mapsto a)}(g(v)) = true$ and $\text{Val}_{B(v \mapsto a)}(s(v)) = w$. Recall that $x \in \text{TC}_B(w)$. But $\text{Pre}(t, B)$ has a node labeled $(s(v), B(v \mapsto a))$ and thus there is a node labeled $(s(v), B(v \mapsto a))$ below $n_0$ in $\text{Pre}(\Pi, A)$. This contradicts the choice of $n_0$.

Finally, suppose that $t$ is a variable $v$. As $\Pi$ is a program and thus has no free variables, $v$ is bound in $\Pi$. Since $\Pi$ does not reuse variables, $v$ is bound exactly once, either by a $\{s(v) : v \in r : g(v)\}$ construction or by a do-forall. Let $r$ be the range of $v$. By Lemma 1, $B$ must be $C(v \mapsto a)$, where $C$ is an expansion of $A_i$ that involves only variables different from $v$ and where $a \in \text{Val}_C(r)$. So $\text{Val}_B(t) = \text{Val}_{C(v \mapsto a)}(v) = a \in \text{Val}_C(r)$. Thus, $x \in \text{TC}(\text{Val}_C(r))$. But $\text{Pre}(\Pi, A)$ includes a copy of $\text{Pre}(r, C)$, whose top node is labeled with $(r, C)$, *below* node $n_0$. This contradicts the minimality of $n_0$.

$\square$

**Corollary 2** *Let $\bar{\Pi}$ be a PTime program $(\Pi, p(n), q(n))$, and let $\rho$ be the run of $\bar{\Pi}$ on some input structure $I$. The number of objects active in $\rho$ is bounded by the number of microsteps plus the number of atoms plus two.*

**Proof** Except for $0, 1$, and atoms, everything active in $\rho$ is activated in $\rho$. The theorem and the observation preceding it immediately give the desired bound. $\square$

## 5.4 The Robustness of $\tilde{\text{C}}$PTime

We have considered two definitions of PTime programs: the official definition by means of active elements, and the counting-microsteps definition. Even though details of the second definition have been skipped, we have shown in the previous two subsections that the two definitions are equivalent in the sense that they give rise to the same notion of $\tilde{\text{C}}$PTime.

There is another natural definition of PTime programs. Fix a program $\Pi$, and call an object $x$ *relevant* to a state $A$ of $\Pi$ if it is active at $A$ or there exists a dynamic function $f$ such that $x \in \text{TC}(\text{Extent}(f))$. Call1 $x$ is *relevant* to a run $\rho$ of $\Pi$ if it is relevant to some state of $\rho$.

A PTime program can be defined as a pair $(\Pi, r(n))$ where $\Pi$ is a program and $r(n)$ is a polynomial that bounds the number of relevant objects in $\Pi$'s runs. If $\rho$ is the run of $\Pi$ on an initial structure $I$, then the run of $(\Pi, r(n))$ on $I$ is the maximal initial segment $\rho_0$ of $\rho$ such (1) that the number of objects relevant to $\rho_0$ is bounded by $r(n)$, and (2) all states of $\rho_0$ are distinct. The second clause is needed to ensure that $\rho_0$ is finite in the case when $\Pi$ loops on $I$.

**Theorem 3** *The active-object and relevant-object definitions give the same notion of $\tilde{\text{C}}$PTime*

**Proof** First, let $(\Pi, r(n))$ be a PTime program with respect to the relevant-object definition, and let $\rho$ be the run of $(\Pi, r(n))$ on an input structure $I$ of size $n$. Clearly, $r(n)$ bounds the number of active objects in $\Pi$'s runs. It suffices to show that the length of $\rho$ is bounded by a polynomial of $n$ that is independent of $I$.

Let $m$ be the number of dynamic names in $\text{Voc}(\Pi)$, and let $f$ range over dynamic functions of $\Pi$. A state in $\rho$ is uniquely determined by the relevant sets $\text{Extent}(f)$. Hence the number of different states in $\rho$ is at most $r(n)^m$.

Second let $(\Pi, q(n), r(n))$ be a PTime program with respect to the active-object definition, and let $\rho$ be the run of $(\Pi, p(n), q(n))$ on an input structure $I$ of size $n$. It suffices to show that the number of objects relevant to $\rho$ is bounded by a polynomial of $n$ that is independent of $I$.

Let $m$ be the number of dynamic functions in $\Pi$ and let $j$ be the maximum of their arities. A relevant object $x$ has one of the following two forms. First, $x$ may be Extent$(f)$ for some dynamic function $f$. There at most $m \cdot p(n)$ relevant objects of that sort. Second, $x$ may be a $k$-tuple of active objects, $k \leq j+1$, or a member of the transitive closure of such a tuple. Obviously, there is a polynomial bound on the number of such relevant objects. $\square$

# 6 Two Fixed-Point Theorems

## 6.1 Definable Set-Theoretic Functions

The ASM programming language allows one to use much of the usual set-theoretic notation. Here are some examples.

**Lemma 2** *Over ASM states, every first-order formula with bounded quantifiers is expressible by a Boolean term.*

**Proof** An easy induction over the given formula. In particular, $(\exists v \in r)\, g(v) \iff 0 \in \{0 : v \in r : g(v)\}$. $\square$

**Lemma 3** *The function*

$$\textit{if } y \textit{ then } x_1 \textit{ else } x_2 = \begin{cases} x_1 & \textit{if } y \neq 0 \\ x_2 & \textit{if } y = 0 \end{cases}$$

*is definable*

**Proof**

$$\text{TheUnique}\Big(\{v : v \in \{x_1, x_2\} : (y \neq 0 \wedge v = x_1) \vee (y = 0 \wedge v = x_2)\}\Big).$$

$\square$

**Lemma 4** *Operations $x \cup y$, $\bigcap x$, $x - y$ are definable.*

**Proof**

$$\begin{aligned} x \cup y &= \bigcup \{x, y\} \\ \bigcap x &= \{v : v \in \bigcup x : (\forall w \in x) v \in w\} \\ x - y &= \{w : w \in x : w \notin y\} \end{aligned}$$

□

The standard Kuratowski definition of ordered pairs is

$$\text{OP}(x, y) = \{\{x\}, \{x, y\}\}.$$

**Lemma 5** *There are definable functions P1 and P2 satisfying the following condition. If $z = OP(x, y)$, then $P1(z) = x$ and $P2(z) = y$.*

**Proof**

$$\begin{aligned} \text{P1}(z) &= \text{TheUnique}\left(\bigcap z\right) \\ \text{P2}(z) &= \left(\text{if } \bigcup z = \bigcap z \text{ then P1}(z) \text{ else TheUnique}(\bigcup z - \bigcap z)\right) \end{aligned}$$

□

We will use the following lemma. Every nonempty transitive set $T$ is a natural model of the vocabulary $\{\in, \emptyset\}$; this model will be also called $T$.

**Lemma 6** *There exists a formula PosInteger(x) in the vocabulary $\{\in, \emptyset\}$ such that, for every transitive set $T$,*

$$T \models PosInteger(x) \iff x \text{ is a positive integer.}$$

**Proof** PosInteger($x$) asserts that $x$ is transitive and either 0 or of the form $z \cup \{z\}$, and the same is true for each $y \in x$. □

## 6.2 First-Order Semantics

The sequel of a given state with respect to a given program can be described in the given state by means of first-order formulas [Glavan and Rosenzweig 1993]. We need here a related result.

**Lemma 7** *For every rule $R$ and every dynamic function name $f$, there is a first-order formula $\text{Update}_{R,f}(\bar{x}, y)$ such that*

$$A \models \text{Update}_{R,f}(\bar{x}, y) \iff (f, \bar{x}, y) \in \text{Den}(R, A)$$

*for all appropriate expanded states $A$.*

The appropriateness of $A$ means that $A$ is appropriate for $R$ and its vocabulary contains the name $f$ which may or may not occur in $R$.

**Proof** Induction on $R$. If $R$ is Skip or an update rule with head name different from $f$, then $\text{Update}_{R,f}$ is any logically false formula. If $R$ is $f(\bar{t}) := t_0$, then $\text{Update}_{R,f}$ is $(\bar{x} = \bar{t} \wedge y = t_0)$. If $R$ is "if $g$ then $R_1$ else $R_2$ endif", then $\text{Update}_{R,f}$ is $(g = \textit{true} \wedge \text{Update}_{R_1,f}) \vee (g = \textit{false} \wedge \text{Update}_{R_2,f})$.

If $R$ is "do-forall $u \in r$, $R_0(u)$ enddo", then $\text{Update}_{R,f}$ is

$$\neg \text{Clash}_R \wedge (\exists u \in r) \text{Update}_{R_0(u),f}$$

where $\text{Clash}_R = \bigvee_{h \in \Upsilon} \text{Clash}_{R,h}$ and $\text{Clash}_{R,h}$ is

$$(\exists u, u' \in r) \exists \bar{z} \exists w \exists w' \left[ \text{Update}_{R^+(u),h}(\bar{z}, w) \wedge \text{Update}_{R^+(u'),h}(\bar{z}, w') \wedge w \neq w' \right]$$

Here $\text{Length}(\bar{z}) = \text{Arity}(h)$. □

## 6.3 Time-Explicit Programs

Call a PTime program $\Pi$ *time-explicit* if every positive integer $i$ is active in all runs of $\Pi$ of length $\geq i$.

**Lemma 8** *Every PTime program can be simulated by a time-explicit PTime program.*

**Proof** Just alter the given program $\Pi$ to

```
do-in-parallel
   Π
   if not(Halt) then
      CT := CT ∪ {CT}
   endif
enddo
```

where CT (an allusion to Current Time) is a fresh nullary dynamic function name (automatically initialized to 0, according to our conventions). □

## 6.4 Fixed-Point Definability

Fix a PTime program $\bar{\Pi} = (\Pi, p(n), q(n))$ where $\Pi$ is time-explicit. Let $I$ range over input structures for $\Pi$. Define $\mathrm{Active}(I)$ to be the collection of active objects in the run of $\bar{\Pi}$ on $I$. It is easy to see that $\mathrm{Active}(I)$ is transitive. We also denote by $\mathrm{Active}(I)$ the structure $(Active(I), \in, \emptyset, \bar{R})$ where $\bar{R}$ stands for all the relations of the input structure $I$.

**Theorem 4 (First Fixed-Point Theorem)** *Let $\langle A_i : i \leq l \rangle$ be the run of $\bar{\Pi}$ on an input structure $I$. Relations*

$$D_f(i, \bar{x}, y) \iff A_i \models f(\bar{x}) = y \neq 0,$$

*where $f$ range over $DynamicVoc(\Pi)$, are uniformly FO+LFP definable in $\mathrm{Active}(I)$.*

The uniformity means that the defining formulas are independent of $I$.

**Proof** Notice that if $i$ is a positive integer then $i - 1 = \bigcup i$. For clarity, we will use $i - 1$ instead of $\bigcup i$ in the situations where $i$ is a positive integer.

Call a first-order formula $\varphi$ *simple* if every atomic subformula of $\varphi$ has the form $f(\bar{x}) = t$ where $\bar{x}$ is a tuple of variables and $t$ is either a variable or *true* or *false*. It is easy to see that every first-order formula whose vocabulary consists of function names is logically equivalent to a simple formula. Without loss of generality, we may assume that the formulas $\mathrm{Update}_{R,f}(\bar{x}, y)$, constructed in the previous subsection, are simple.

By simultaneous recursion, we define relations $D_f$, where $f$ ranges over the dynamic function names of $\Pi$:

$$D_f(i, \bar{x}, y) \iff \mathrm{PosInteger}(i) \wedge \\ \left[ \left( D_f(i-1, \bar{x}, y) \wedge \neg(\exists z \neq y) U_f(i-1, \bar{x}, z) \right) \vee U_f(i-1, \bar{x}, y) \right]$$

Here $U_f(j, \bar{x}, y)$ is the formula $\mathrm{Update}_{\Pi, f}(\bar{x}, y)$ where each atomic subformula $h(\bar{u}) = t$ is replaced with $D_h(j, \bar{u}, t)$. $\square$

The First Fixed-Point Theorem remains true if the computation time of the given program is bounded by any other function (not necessarily a polynomial) or is not bounded at all. Also, we get the same definability in any transitive $T$ that includes $\mathrm{Active}(I)$.

**Theorem 5 (Second Fixed-Point Theorem)** *Restrict attention to input structures $I$ such that $\bar{\Pi}$ halts on $I$. Then the set $Active(I)$ is uniformly FO+LFP definable in $HF(BaseSet(I))$.*

**Proof** The desired FO+LFP formula $\varphi(x)$ asserts that $x$ is an atom, or $x \in \{0,1\}$, or the following is true for some dynamic function $f$ where $k = \text{Arity}(f)$.

$$(\exists i, v_0, \ldots, v_k)\Big[D_f(i, v_0, \ldots, v_k) \wedge (v_0 = x \vee \cdots \vee v_k = x)\Big]$$

□

# 7 On the Extent of $\tilde{\text{C}}$PTime

We show, in particular, that PTime abstract state machines are more powerful than the PTime relational machines of Abiteboul–Vianu.

**Theorem 6** *For every PTime relational machine $\Xi$, there exists a PTime ASM program $\Pi$ that accepts all input structures accepted by $\Xi$ and rejects all input structures rejected by $\Xi$.*

**Proof** If $\Xi$ has $m$ instructions, then the desired program $\Pi$ is a do-in-parallel rule with $m$ components. Each component simulates one instruction of $\Xi$. □

Let $\Upsilon_0$ be an input vocabulary and let $\Upsilon_1$ be the extension of $\Upsilon_0$ with a unary predicate $U$. If $A$ is a structure with basic predicate $U$, let $A_0$ be the $\Upsilon_0$-reduct of the substructure of $A$ with universe $U$.

**Theorem 7** *Let $j$ be a positive integer and let $K$ be the class of finite $\Upsilon_1$-structures $A$ where $|A_0|! \leq |A|^j$. For every PTime class $K_0$ of $\Upsilon_0$-structures, there exists an ASM program $\Pi$ with input vocabulary $\Upsilon_1$ that accepts all $K$-structures $A$ with $A_0$ in $K_0$ and rejects all $K$-structures $A$ with $A_0$ outside of $K_0$.*

**Proof** Suppose that $A \in K$ and let $n = |A_0|$. There exists a Turing machine $T$ that, given an ordered version of $A_0$, computes $F(A|U)$. The desired ASM program $\Pi$ constructs all $n!$ orderings of $A_0$ in parallel and then runs in parallel $n!$ simulations of $T$. □

In particular, we have the case when $\Upsilon_0$ is empty, $j = 1$ and $K_0$ is the class of naked sets of even cardinality.

**Lemma 9** *Let $K$ be the class of all $\{U\}$ structures $A$ with $|A_0|! \leq |A|$. There exists no relational machine that accepts all $K$-structures $A$ with $|A_0|$ even and rejects all $K$-structures $A$ with $|A_0|$ odd.*

**Proof** Relational machines recognize only classes definable in FO+LFP [Abiteboul and Vianu 1991]. It is easy to see that $K$ admits quantifier elimination and that every FO+LFP formula is equivalent to a quantifier-free first-order formula. It follows that parity is not FO+LFP definable. $\square$

Thus PTime abstract state machines are more powerful than PTime relational machines.

Remark. In the case when $\Upsilon_0$ is empty, $j = 1$ and $K_0$ is the class of naked structures of even cardinality, the theorem above can be strengthened by replacing the restriction $|A_0|! \leq |A|^j$ with more liberal restriction $2^{|A_0|} \leq |A|^j$. The idea is to compute the set of pairs of 2-element subsets of $U$, then extend it with the set of all 4-element subsets of $U$, then extend the result with the set of all 6-element subsets of $U$, and so on. When this computation converges, check if the result contains $U$.

Remark. Theorem 7 and its proof can be extended to cover the situation where $U$ is not merely a subset of the input structure but rather a set that can be produced in polynomial time by an ASM. For example, if the intput structures are graphs $G$ then $U$ might be the commutator subgroup $G'$ of the quotient $G/G'$.

## 8 The Support Theorem

The goal of this and the next sections is to show that the parity of a naked set is not ČPTime computable. Thus the inclusion of ČPTime in PTime (see Theorem 1) is proper; "choiceless" is a real restriction. The present section is devoted to establishing a limitation on the sets that can be activated by a ČPTime computation over a naked set. This limitation is used in the next section to prove the negative result about parity. The same method will also yield other negative results.

Consider a PTime program $\Pi$ and let $I$ be an input structure for $\Pi$. The recipe $\theta(x) = \{\theta(y) : y \in x\}$ extends any automorphism $\theta$ of $I$ to an automorphism of the whole initial state $\text{State}(I)$ generated by $I$. It is easy to see that every automorphism of $\text{State}(I)$ can be obtained this way. Indeed, an automorphism $\theta'$ of $\text{State}(I)$ coincides, on $I$, with some automorphism $\theta$ of $I$; by induction on $\text{Rank}(x)$, check that $\theta'(x) = \theta(x)$ for all $x \in \text{State}(I)$.

**Definition 3** *A set $X$ of atoms of $I$ is a support of an object $y \in \text{State}(I)$ if every automorphism of $I$ that pointwise fixes $X$ fixes $y$ as well.* $\square$

Let Active($I$) be the set of active objects in the run of $\Pi$ on $I$. It is easy to see that Active($I$) is transitive and closed under automorphisms of State($I$). Let Active$^+$($I$) be the substructure of State($I$) with base set Active($I$).

**Theorem 8 (Support Theorem)** *Assume that the input vocabulary of $\Pi$ is empty. There exists a number $k$ such that, for all sufficiently large $I$, every object in Active($I$) has a support of cardinality $\leq k$.*

To avoid interruption of the natural flow of the proof, we start with a known (except we do not know a proper reference) combinatorial lemma which will be used later in the proof. Recall that a $\Delta$-*system* is a collection $K$ of sets such that $X \cap Y$ is the same set for all $X \neq Y$ in $K$.

**Lemma 10** *Any family $F$ of $\geq l! p^{l+1}$ sets, each of size $\leq l$, includes a $\Delta$-system of $p$ sets.*

**Proof** Induction on $l$. If $l = 0$, then $F$ itself is a $\Delta$-system. Assume that $l > 0$ and the results holds for $l - 1$.

Case 1: There exists a point $x$ that belongs to $\geq (l-1)! p^l$ sets in $F$, say sets $X_i$, $i \in I$. Apply the induction hypothesis to the family $\{X_i - \{x\} : i \in I\}$, to extract a $\Delta$-system of $p$ sets $\{X_i - \{x\} : i \in J\}$. The family $\{X_i : i \in J\}$ is the desired $\Delta$-system.

Case 2: Each point belongs to $< (l-1)! p^l$ sets in $F$. In this case, we find $p$ pairwise disjoint members of $F$; they form the desired $\Delta$-system. Notice that each member of $F$ intersects $< l(l-1)! p^l = l! p^l$ other members and that $\text{Card}(F)/(l! p^l) \geq p$. Pick a member $X_1$ arbitrarily, and then eliminate those members that meet $X_1$. Pick a member $X_2$ among the remaining members arbitrarily, and then eliminate those members that meet $X_2$. And so on. □

Let $I$ be a naked set and let $\mathcal{A} = \text{Active}(I)$.

**Lemma 11** *If $X_1, X_2$ support $y$ and $X_1 \cup X_2 \neq I$, then $X_1 \cap X_2$ supports $y$ as well.*

**Proof** Suppose that $X_1, X_2$ support $y$. Fix an atom $a \in I - (X_1 \cup X_2)$. Let $b$ range over $I - (X_1 \cap X_2)$ and $\pi_b$ be the transposition of atoms that takes $a$ to $b$. For each $b$, either $b \notin X_1$ or $b \notin X_2$. In the first case $\pi_b$ pointwise fixes $X_1$, and in the second it pointwise fixes $X_2$. In either case, it fixes $y$. It is easy to see that the transpositions $\pi_b$ generate all permutations of atoms which pointwise fix $X_1 \cap X_2$. Hence the automorphisms induced by permutations $\pi_b$ generate all automorphisms of $\mathcal{A}$ that fix $X_1 \cap X_2$. Hence every such automorphism fixes $y$. □

Let $n = \mathrm{Card}(I)$. Lemma 11 justifies the following definition. If object $y$ has a support $X$ with $|X| < n/2$, then the set

$$\mathrm{Supp}(y) = \bigcap \{X : X \text{ supports } y \text{ and } |X| < n/2\}.$$

is the smallest support of $X$.

Since $\Pi$ is PTime, there exists a bound $n^k$ on $\mathrm{Card}(\mathcal{A})$. Fix such a $k$ and assume that $n$ is so large that $\binom{n}{k+1} > n^k$.

**Lemma 12** *If $x \in \mathcal{A}$ has a support of size $< n/2$, then $|Supp(x)| \leq k$.*

**Proof** Suppose that $x$ has a support of size $< n/2$ and $s = |\mathrm{Supp}(x)|$. If an automorphism $\theta$ moves $x$ to some $y$, then it moves $\mathrm{Supp}(x)$ to $\mathrm{Supp}(y)$. If $s > k$, we have

$$\begin{aligned} n^k & \geq \mathrm{Card}(\mathcal{A}) \geq \mathrm{Card}\{\theta(x) : \theta \in \mathrm{Aut}(\mathcal{A})\} \\ & \geq \mathrm{Card}\{\theta(\mathrm{Supp}(x)) : \theta \in \mathrm{Aut}(\mathcal{A})\} \\ & = \binom{n}{s} \geq \binom{n}{k+1} > n^k \end{aligned}$$

□

In order to prove the theorem, it suffices to prove the following lemma.

**Lemma 13** *If $n$ is sufficiently large, then every member of $\mathcal{A}$ has a support of size $< n/2$.*

**Proof** Toward a contradiction, assume that the lemma fails and let $x$ be an object of minimal rank without support of size $< n/2$. Clearly, $x$ is a set and each member of $x$ has a support of size $< n/2$. Let $m = \lfloor n/(4k) \rfloor$.

**Claim 1** *There exists a sequence $\langle (\theta_j, y_j, z_j, Y_j, Z_j) : 1 \leq j \leq m \rangle$ such that every initial segment $\langle (\theta_i, y_i, z_i, Y_i, Z_i) : 1 \leq i \leq j \rangle$ satisfies the following conditions:*

- *$\theta_j$ is an automorphism of $\mathcal{A}$, and $y_j, z_j$ are objects in $\mathcal{A}$, and $Y_j = Supp(y_j)$, $Z_j = Supp(z_j)$.*
- *$y_j \in x$, $z_j \notin x$.*
- *$\theta_j$ fixes $Y_i \cup Z_i$ pointwise for all $i < j$, and $\theta_j(y_j) = z_j$, and $\theta_j$ maps $Y_j$ onto $Z_j$.*

**Proof** of the claim. We construct the tuples by induction on $j$. Suppose that a sequence $\langle (\theta_i, y_i, z_i, Y_i, Z_i) : 1 \leq i < j \rangle$, satisfying all the conditions, has been constructed. By the minimality of $x$, each $y_i$ has a support of size $< n/2$. Since $z_i$ is an automorphic image of $y_i$, the same applies to $z_i$. By the previous lemma, $|Y_i|, |Z_i| \leq k$.

Let $X_j = \bigcup_{i<j}(Y_i \cup Z_i)$. We have $|X_j| \leq (j-1) \cdot 2k < (n/4k) \cdot 2k = n/2$. If every automorhism $\theta$ that pointwise fixes $\bigcup_{i<j}(Y_i \cup Z_i)$ fixes $x$ as well, then $x$ has a support of size $< n/2$ and we have the desired contradiction. Suppose that there exists an automorphism $\theta$ that pointwise fixes $X_j$ but moves $x$. It follows that there exists $y \in x$ such that the element $z = \theta(y)$ does not belong to $x$. (Otherwise $\theta(x) = \theta\{y : y \in x\} = \{\theta(y) : y \in x\} = x$.) Since $\theta(y) = z$, $\theta$ maps $\mathrm{Supp}(y)$ onto $\mathrm{Supp}(z)$. Choose, $\theta_j = \theta, y_j = y$ and $z_j = z$. □

Let $p$ be the largest integer with $(2k)! p^{2k+1} \leq m$. As $n$ grows, both $m$ and $p$ grow (but $k$ is fixed). For large enough $n$, we have

$$2^{p-1} > \left[\left((2k)!(p+1)^{2k+1} \cdot 4k\right)^k\right] \geq \left[\left((m+1) \cdot 4k\right)^k\right] > n^k$$

Assume that $n$ is sufficiently large, so that $2^{p-1} > n^k$.

**Claim 2** *There exists a sequence $\langle (\theta_j, y_j, z_j, Y_j, Z_j) : 1 \leq j \leq p \rangle$ such that*

- *every initial segment $\langle (\theta_i, y_i, z_i, Y_i, Z_i) : 1 \leq i \leq j \rangle$ of the given sequence satisfies the three conditions of Claim 1, and*
- *The sets $Y_i \cup Z_i$ form a $\Delta$-system.*

**Proof** of the claim. Let $\langle (\theta_j, y_j, z_j, Y_j, Z_j) : 1 \leq j \leq m \rangle$ be as in Claim 1. By the induction hypothesis, eachy $Y_i$ is of cardinality $\leq k$. Since $Z_i$ is an automorphic image of $Y_i$, the same applies to $Z_i$. Thus $Y_i \cup Z_i \leq 2k$. Now apply Lemma 10. □

Fix a sequence $\langle (\theta_j, y_j, z_j, Y_j, Z_j) : 1 \leq j \leq p \rangle$ as in Claim 2, and let $X_0 = (Y_i \cup Z_i) \cap (Y_j \cup Z_j)$ for all $i \neq j$ in $[1, \ldots, p]$. Let $U$ be the integer interval $[2..p]$. If $i \in U$, then $\theta_i$ pointwise fixes $Y_1 \cup Z_1$ and therefore pointwise fixes $X_0$.

**Claim 3** *For each $V \subseteq U$, there exists an automorphism $\theta_V$ such that*

- *if $i \in V$, then $z_i = \theta_V(z_i)$.*
- *if $i \in U - V$, then $z_i = \theta_V(y_i)$.*

**Proof** of the claim. Construct a permutation $\pi(a)$ of atoms as follows. If $a \in X_0$ then $\pi(a) = a$. If $a \in Y_i \cup Z_i$ for some $i \in V$, then $\pi(a) = a$. If $a \in Y_i$ for some $i \in U - V$ but $a \notin X_0$, then $\pi(a) = \theta_i(a)$, so that $\pi$ maps $Y_i$ onto $Z_i$. We do not care how $\pi$ behaves on the remaining atoms. The desired $\theta_V$ is the automorphism induced by $\pi$. □

Let the automorphisms $\theta_V$ be as in Claim 12.

**Claim 4** *If $V, W$ are different subsets of $U$, then $\theta_V(x) \neq \theta_W(x)$.*

**Proof** of the claim. Suppose that $V$ and $W$ are distinct. Without loss of generality, $V - W \neq \emptyset$. Pick some $i \in V - W$. We show that $z_i \in \theta_W(x) - \theta_V(x)$.

Since $z_i \notin x$, $\theta_V(z_i) \notin \theta_V(x)$. Since $i \in V$, $z_i = \theta_V(z_i) \notin \theta_V(x)$.

Since $y_i \in x$, $\theta_W(y_i) \in \theta_W(x)$. Since $i \in U - W$, $z_i = \theta_W(y_i) \in \theta_W(x)$. □

By Claim 3, there are $2^{p-1}$ different automorphic images of $x$. Recall that $n$ is sufficiently large, so that $2^{p-1} > n^k$. Hence $\text{Card}(\mathcal{A}) \geq 2^{p-1} > n^k \geq \text{Card}(\mathcal{A})$. This gives the desired contradiction. The Support Theorem is proved. □

Define a *colored set* to be an input structure with only unary relations, called *colors*, which partition the base set, so that every atom (that is every element of the base set) belongs to exactly one color.

**Corollary 3** *Assume that input structures for the PTime program $\Pi$ are colored sets. There exists a number $k$, depending only on $\Pi$, satisfying the following condition. If every color in a given input structure $I$ is sufficiently large, then every object in $\text{Active}_\Pi(I)$ has a support of cardinality $\leq k$.*

**Proof** The proof is similar to the proof of the theorem. We indicate the more important changes. In Lemma 11, require that $I - (X_1 \cup X_2)$ contains at least one atom of every color.

In the definition of $\text{Supp}(y)$, replace $|X| < n/2$ with the requirement

$$|X \cap C| < |C|/2 \quad \text{for all colors } C \qquad (*)$$

Accordingly, in Lemma 12 and Lemma 13, instead of supports of size $< n/2$, speak about supports satisfying $(*)$. In the definition of $m$, replace $n$ with the minimum of the color sizes. The rest of the proof remains valid. □

Finally, let us note that, over some input structures, a PTime program can activate sets with no bounded support, so that the minimal support size depends on the input structure $I$, not only on the PTime program.

**Example**  Let $I$ be the disjoint union of (i) a vector space $V$ over the two-element field and (ii) a disjoint set $S$ of size $\geq 2^{|V|}$. Let the program do the following with a dynamic nullary function $Q$. Initialize $Q$ to $\{\{\bar{0}\}\}$ where $\bar{0}$ is the zero of $V$. Thereafter, for each $q \in Q$ and each $v \in V - q$, put into $Q$ the subspace generated by $q \cup \{v\}$, *except* if this would make $V \in Q$, in which case halt and accept. On the first step, the program activates all one-dimensional subspaces of $V$; on the second, all two-dimensional subspaces; on the third, all three-dimensional subspaces, and so on. The length of the run equals the dimension of $V$. The number of active objects in the run is

$$|V| + |S| + |\text{Subspaces of}(V)| \leq |V| + 2|S| < 2|I|.$$

Thus a PTime version of $\Pi$ accepts $I$. Notice that every hyperplane $H$ of $V$ is activated. But $H$ has no support smaller than $\dim(V) - 1$. □

## 9  The Equivalence Theorem

Fix an input vocabulary $\Upsilon_0$ and let $I, J$ denote input structures of vocabulary $\Upsilon_0$. Recall that every automorphism of $I$ naturally extends to an automorphism of $\mathrm{HF}(I)$ and that a subset $X$ of $\mathrm{BaseSet}(I)$ supports an object $y \in \mathrm{HF}(I)$ if every automorphism of $I$ that pointwise fixes $X$ also fixes $y$. Given a positive integer $k$, call an object $y \in \mathrm{HF}(I)$ *k-symmetric* if every $z \in \mathrm{TC}(y)$ has a support of size $\leq k$. This terminology makes sense in the case of interest to us when many permutations of $I$ fix a $k$-symmetric object $y$. Notice that every atom has a support of size one and thus is $k$-symmetric. Let $\bar{I}_k$ denote the collection of $k$-symmetric objects in $\mathrm{HF}(I)$ as well as the corresponding structure of vocabulary $\Upsilon_0 \cup \{\in, \emptyset\}$.

We are interested in a special case when $\Upsilon_0$ is impty and thus $I, J$ are naked sets.

**Theorem 9 (Equivalence Theorem)** *Fix positive integers $k$ and $m$. If naked sets $I, J$ are sufficiently large, then structures $\bar{I}_k$ and $\bar{J}_k$ are $L^m_{\infty,\omega}$-equivalent.*

The theorem is proved in the rest of this section. We drop the subscript $k$ and abbreviate "$k$-symmetric" to "symmetric". Without loss of generality, $m \geq 3$. We assume that the naked sets $I, J$ have size $\geq km$ and construct a winning strategy for the Duplicator in $\mathrm{Game}_m(\bar{I}, \bar{J})$. The idea is to represent every relevant object $x$ as a combination of a form and matter. The form of an object $x$ reflects a definition of $x$ independent from the underlying sets of atoms. The matter of $x$ is an ordered support of $x$.

## 9.1 Matter

**Molecules** A *molecule* over a naked set $I$ is an injective map $\sigma : k \longrightarrow I$. In other words, a molecule is a sequence of $k$ distinct atoms.

**The Configuration of a Sequence of Molecules** Consider a naked set $I$. The *configuration* $C(\bar{\sigma})$ of a finite sequence $\bar{\sigma} = (\sigma_0, \ldots, \sigma_{l-1})$ of molecules over $I$ is the equivalence relation on $l \times k$ given by

$$(i,p)C(\bar{\sigma})(j,q) \iff \sigma_i(p) = \sigma_j(q).$$

The configuration describes how the ranges of the molecules overlap. By the injectivity, $(i,p)C(\bar{\sigma})(i,q) \iff p = q$. Notice that $C(\bar{\sigma})$ is uniquely determined by the configurations $C(\sigma_i, \sigma_j)$.

**Abstract Configurations** Let $l$ be a natural number. An *$l$-ary configuration* is an equivalence relation $E$ on $l \times k$ such that $(i,p)E(i,q) \iff p = q$. Every $C(\sigma_0, \ldots, \sigma_{l-1})$ is an $l$-ary configuration, and every $l$-ary configuration can be realized in this way. To prove the latter, assign a different atom $[i,p]$ to every equivalence class $(i,p)_E$ of $E$. Then construct $\sigma_i(p) = [i,p]$.

**Lemma 14** *Suppose that*

  *$l \in m$, and $\sigma_0, \ldots, \sigma_l$ are molecules over $I$, and $\tau_1, \ldots, \tau_l$ are molecules over $J$;*

  *$Q = C(\sigma_0, \ldots, \sigma_l)$ and $Q' = C(\sigma_1, \ldots, \sigma_l) = C(\tau_1, \ldots, \tau_l)$.*

*There exists a molecule $\tau_0$ over $J$ with $C(\tau_0, \ldots, \tau_l) = Q$.*

**Proof** Define the desired $\tau_0$ by setting

$$\tau'_0(p) = \begin{cases} \tau_1(q_1) & \text{if } (0,p)Q(1,q_1) \\ \tau_2(q_2) & \text{if } (0,p)Q(2,q_2) \\ \ldots \\ \tau_l(q_l) & \text{if } (0,p)Q(l,q_l) \end{cases}$$

and then extending $\tau'_0$ to a full molecule $\tau_0$ by using distinct values in $J - \bigcup_{i=1}^{l} \text{Range}(\tau_i)$. It is obvious that $C(\tau_0, \ldots, \tau_l) = Q$ provided that $\tau_0$ is well-defined. Recall that $\text{Card}(J) \geq km$ and thus there exist enough distinct values to extend $\tau'_0$ to $\tau_0$. In the rest of the proof, we check that $\tau'_0$ is well-defined.

First, we check that each $\tau'_0(p)$ is defined uniquely. Let $1 \leq i, j \leq l$ and suppose $(0,p)Q(i,q)$ and $(0,p)Q(j,s)$. Then $(i,q)Q(j,s)$, $(i,q)Q'(j,s)$, and therefore $\tau_i(q) = \tau_j(s)$.

Second, we check that $\tau'_0$ is injective. Let $1 \leq i, j \leq l$ and suppose that $\tau'_0(p) = \tau'_0(p')$ where $\tau'_0(p) = \tau_i(q)$ and $\tau'_0(p') = \tau_j(s)$. By the definition of $\tau'_0$, we have $(0, p)Q(i, q)$ and $(0, p')Q(j, s)$. Further, $\tau_i(q) = \tau'_0(p) = \tau'_0(p') = \tau_j(s)$, and hence $(i, q)Q'(j, s)$ and therefore $(i, q)Q(j, s)$. Putting this together, we have

$$(0, p) \ Q \ (i, q) \ Q \ (j, s) \ Q \ (0, p')$$

which implies $p = p'$. □

## 9.2 Forms

**The Definition of Forms** Fix a list $c_0, c_1, \ldots, c_{k-1}$ of new symbols. The set of *forms* is the smallest set containing the symbols $c_p$ and containing every finite set of pairs $(\varphi, E)$ where $\varphi$ is a form and $E$ is a binary configuration. If $\varphi = c_p$ then $\mathrm{Rank}(\varphi) = 0$; otherwise

$$\mathrm{Rank}(\varphi) = 1 + \max\{\mathrm{Rank}(\psi) : \text{ some } (\psi, E) \in \varphi\}.$$

**Denotations** A form $\varphi$ and a molecule $\sigma$ over a naked set $I$ uniquely define an object $\varphi *_I \sigma \in \mathrm{HF}(I)$. The subscript may be omitted in $*_I$ if the naked set is clear from the context. The definition is given by induction on $\varphi$.

- $c_p * \sigma = \sigma(p)$.

- If $\varphi$ is a set, then $\varphi * \sigma = \{\psi * \tau : (\psi, C(\tau, \sigma)) \in \varphi\}$.

**Permutations** Any permutation $\pi$ of $I$ extends to an automorphism, also called $\pi$, of the structure $(\mathrm{HF}(I), \in)$ by means of the following rule: $\pi(x) = \{\pi(y) : y \in x\}$. It is easy to see that every automorphism of $\mathrm{HF}(I)$ is obtained this way.

**Lemma 15** *If $\pi$ is a permutation of $I$, then $\pi(\varphi * \sigma) = \varphi * \pi\sigma$.*

**Proof** by induction on $\varphi$. $\pi(c_p * \sigma) = \pi\sigma(p) = (\pi\sigma)(p) = c_p * (\pi\sigma)$. If $\varphi$ is a set, then

$$\begin{aligned}
\pi(\varphi * \sigma) &= \pi\{\psi * \tau : (\psi, C(\tau, \sigma)) \in \varphi\} = \{\pi(\psi * \tau) : (\psi, C(\tau, \sigma)) \in \varphi\} \\
&= \{\psi * \pi\tau : (\psi, C(\tau, \sigma)) \in \varphi\} \quad \text{(by induction hypothesis)} \\
&= \{\psi * \rho : (\psi, C(\pi^{-1}\rho, \sigma)) \in \varphi\} \quad (\rho = \pi\tau) \\
&= \{\psi * \rho : (\psi, C(\rho, \pi\sigma)) \in \varphi\} \quad \text{(see below)} \\
&= \varphi * \pi\sigma
\end{aligned}$$

It remains to verify that $C(\pi^{-1}\rho, \sigma) = C(\rho, \pi\sigma)$:

$$(0,p)C(\pi^{-1}\rho,\sigma)(1,q) \iff (\pi^{-1}\rho)(p) = \sigma(q) \iff$$
$$\rho(p) = (\pi\sigma)(q) \iff (0,p)C(\rho,\pi\sigma)(1,q)$$

□

**Corollary 4** *Every $\varphi *_I \sigma$ is symmetric and thus belongs to $\bar{I}$.*

**Proof** Indeed, every $\varphi *_I \sigma$ is supported by $\text{Range}(\sigma)$ and thus has a support of size $\leq k$. The same conclusion applies to the members of $\varphi *_I \sigma$. □

Recall that $I$ is a naked set of cardinality $\geq km$.

**Lemma 16** *Every symmetric object $x$ over $I$ is equal to $\varphi *_I \sigma$ for some form $\varphi$ and some molecule $\sigma$ over $I$.*

**Proof** Any atom $x$ equals $c_0 * \sigma$ where $\sigma$ is an arbitrary molecule with $\sigma(0) = x$. Proceeding inductively, suppose that $x$ is a symmetric set with elements $y = \psi_y * \tau_y$. Since $x$ is symmetric, there is a molecule $\sigma$ whose range supports $x$. We will prove that $x = \varphi * \sigma$ where $\varphi = \{(\psi_y, C(\tau_y, \sigma)) : y \in x\}$. One inclusion is easy. Suppose that $y \in x$. By the definition of $*_I$, $\varphi * \sigma = \{\psi * \tau : (\psi, C(\tau, \sigma)) \in \varphi\}$. By the definition of $\varphi$, $(\psi_y * C(\tau_y, \sigma)) \in \varphi$. Hence $y = \psi_y * \tau_y \in \varphi * \sigma$.

For the difficult direction, consider any $z \in \varphi * \sigma$. By the definition of $*_I$, $z$ is a composition $\psi * \rho$ such that $(\psi, C(\rho, \sigma)) \in \varphi$. By the definition of $\varphi$, there exists a $y$ such that $\psi = \psi_y$ and $C(\rho, \sigma) = C(\tau_y, \sigma)$. The latter equality does not imply that $\rho = \tau_y$ and we are not going to prove that $z = y$. Instead we construct an automorphism $\pi$ of $\bar{I}$ that pointwise fixes $\sigma$ and moves $y$ to $z$. Since $\sigma$ supports $x$, $\pi$ fixes $x$; hence $z \in x$. It remains to construct such $\pi$.

We want that $\pi\tau_y = \rho$ and that $\pi$ pointwise fixes $\text{Range}(\sigma)$. To this end, define a function $\pi_0 : \text{Range}(\tau_y) \cup \text{Range}(\sigma) \longrightarrow I$ by

$$\pi_0(a) = \begin{cases} \rho(p) & \text{if } a = \tau_y(p); \\ a & \text{if } a = \sigma(q) \end{cases}$$

Even though the two cases are not mutually exclusive, $\pi_0$ is well-defined. Indeed,

$$\tau_y(p) = \sigma(q) \Rightarrow (0,p)C(\tau_y,\sigma)(1,q) \Rightarrow (0,p)C(\rho,\sigma)(1,q) \Rightarrow \rho(p) = \sigma(q).$$

Furthermore, $\pi_0$ is injective. Indeed, assume that $\pi_0(a) = \pi_0(b)$. If $\pi_0(a) = \rho(p_1), \pi_0(b) = \rho(p_2)$ then $p_1 = p_2$ (because $\rho$ is injective) and therefore $a = \tau_y(p_1) = \tau_y(p_2) = b$. In case $a = \sigma(q_1), b = \sigma(q_2)$, we have $a = \pi_0(a) = \pi_0(b) = b$. Finally suppose that $a = \tau_y(p), b = \sigma(q)$. Then

$$\pi_0(a) = \pi_0(b) \Rightarrow \rho(p) = \sigma(q) \Rightarrow (0,p)C(\rho,\sigma)(1,q) \Rightarrow$$
$$(0,p)C(\tau_y,\sigma)(1,q) \Rightarrow \tau_y(p) = \sigma(q) \Rightarrow a = b$$

Thus, function $\pi_0$ is one-to-one. Extend it to a permutation $\pi$ over $I$ in an arbitrary way. Since $\pi$ extends $\pi_0$, it pointwise fixes Range$(\sigma)$ and $\pi\tau_y = \rho$. In the standard way, $\pi$ extends to an automorphism of $\bar{I}$ which will be denoted $\pi$ as well. Since $\sigma$ supports $x$ (by the choice of $\sigma$), $\pi(x) = x$. By Lemma 15,

$$\pi(y) = \pi(\psi_y * \tau_y) = \psi_y * (\pi\tau_y) = \psi_y * \rho = z.$$

□

## 9.3 The In and Eq Relations

**Lemma 17** *There are ternary relations Eq and In such that, in every $\bar{I}$,*

$$\psi * \tau \in \varphi * \sigma \iff In(\psi, \varphi, C(\tau, \sigma)) \qquad (1)$$
$$\psi * \tau = \varphi * \sigma \iff Eq(\psi, \varphi, C(\tau, \sigma)) \qquad (2)$$

*for all forms $\varphi, \psi$ and all molecules $\sigma, \tau$.*

The crucial points here are that $\sigma$ and $\tau$ are involved in In and Eq only via their configurations and that In and Eq don't depend on $I$.

**Proof** We define In$(\psi, \varphi, E)$ and Eq$(\psi, \varphi, E)$ by recursion on Rank$(\psi)$ + Rank$(\varphi)$.

$$\begin{aligned}
\text{In}(\psi, \varphi, E) \iff & \;\varphi \text{ is a set and } (\exists \text{ form } \chi) \\
& (\exists \text{ ternary configuration } Q \text{ with } Q_{12} = E) \\
& \big[(\chi, Q_{02}) \in \varphi) \text{ and } \text{Eq}(\psi, \chi, Q_{10})\big]
\end{aligned}$$

$$\begin{aligned}
\text{Eq}(\psi, \varphi, E) \iff & \;\text{either } (\exists p, q \in k)\big[\psi = c_p \wedge \varphi = c_q \wedge (0,p)E(1,q)\big], \\
& \text{or } \varphi, \psi \text{ are sets and } (\forall \text{ form } \chi) \\
& (\forall \text{ ternary configuration } Q \text{ with } Q_{12} = E) \\
& \big[\text{if } (\chi, Q_{02}) \in \varphi \text{ then } \text{In}(\chi, \psi, Q_{01}), \text{ and} \\
& \text{if } (\chi, Q_{01}) \in \psi \text{ then } \text{In}(\chi, \varphi, Q_{02})\big]
\end{aligned}$$

Proof of (1). If $\varphi$ is a symbol $c_p$, then $\varphi * \sigma$ is an atom, so the left side of (1) is false. So is the right side, by the definition of In. Thus, we may assume from now on that $\varphi$ is a set.

Suppose first that $\psi * \tau \in \varphi * \sigma$. By the definition of $*_I$, $\psi * \tau = \chi * \rho$ for some $\chi, \rho$ with $(\chi, C(\rho, \sigma)) \in \varphi$. By the induction hypothesis, $\mathrm{Eq}(\psi, \chi, C(\tau, \rho))$. We check that this $\chi$ and the ternary configuration $Q = C(\rho, \tau, \sigma)$ witness $\mathrm{In}(\psi, \varphi, C(\tau, \sigma))$. Indeed, $Q_{12} = C(\tau, \sigma)$, $(\chi, Q_{02}) = (\chi, C(\rho, \sigma)) \in \varphi$, and $\mathrm{Eq}(\psi, \chi, Q_{10})$ is $\mathrm{Eq}(\psi, \chi, C(\tau, \rho))$.

Conversely, suppose that $\mathrm{In}(\psi, \varphi, C(\tau, \sigma))$ is witnessed by $\chi$ and $Q$. By Lemma 14, there exists $\rho$ such that $Q = C(\rho, \tau, \sigma)$. We have:

$(\chi, C(\rho, \sigma)) = (\chi, Q_{02}) \in \varphi$, so that $\chi * \rho \in \varphi * \sigma$; and

$\mathrm{Eq}(\psi, \chi, Q_{10})$ holds, that is $\mathrm{Eq}(\psi, \chi, C(\tau, \rho))$ holds.

By the induction hypothesis, $\psi * \tau = \chi * \rho \in \varphi * \sigma$. Part (1) is proved.

Proof of (2). Both sides of (2) are false if one of $\psi, \varphi$ is a symbol $c_p$ while the other is a set. If $\varphi = c_q, \psi = c_p$, then

$$\psi * \tau = \varphi * \sigma \iff \tau(p) = \sigma(q) \iff (0,p)C(\tau,\sigma)(1,q)$$
$$\iff \mathrm{Eq}(\psi, \varphi, C(\tau, \sigma))$$

So we may assume from now on that both $\psi$ and $\varphi$ are sets.

Suppose first that $\psi * \tau = \varphi * \sigma$. Let $\chi$ be any form and $Q$ be any ternary configuration with $Q_{12} = C(\tau, \sigma)$. We must prove

$$(\chi, Q_{02}) \in \varphi \implies \mathrm{In}(\chi, \psi, Q_{01})$$
$$(\chi, Q_{01}) \in \varphi \implies \mathrm{In}(\chi, \psi, Q_{02})$$

By symmetry, it suffices to prove only the first of these two implications. So assume $(\chi, Q_{02}) \in \varphi$. By Lemma 14, there exists $\rho$ such that $C(\rho, \tau, \sigma) = Q$. Then $(\chi, C(\rho, \sigma)) = (\chi, Q_{02}) \in \varphi$, so $\chi * \rho \in \varphi * \sigma = \psi * \tau$. By the induction hypothesis, $\mathrm{In}(\chi, \psi, C(\rho, \tau))$, that is $\mathrm{In}(\chi, \psi, Q_{01})$, as required.

Conversely, suppose that $\mathrm{Eq}(\psi, \varphi, C(\tau, \sigma))$ holds. By symmetry, it suffices to prove only that $\psi * \tau \subseteq \varphi * \sigma$. Let $\chi * \rho$, with $(\chi, C(\rho, \tau)) \in \psi$, be an arbitrary element of $\psi * \tau$. Apply the definition of $\mathrm{Eq}(\psi, \varphi, C(\tau, \sigma))$ with this $\chi$ and with $Q = C(\rho, \tau, \sigma)$, which satisfies $Q_{12} = C(\tau, \sigma)$. Since $Q_{01} = (\chi, C(\rho, \tau)) \in \psi$, we have $\mathrm{In}(\chi, \varphi, Q_{02})$, that is $\mathrm{In}(\chi, \varphi, C(\rho, \sigma))$. By the induction hypothesis, $\chi * \rho \in \varphi * \sigma$, as required. □

## 9.4 The Winning Strategy

Now we are ready to construct a winning strategy for the Duplicator in the $\text{Game}_m(\bar{I}, \bar{J})$. The strategy is to ensure that, after every step, there exist forms $\varphi_i$, $I$-molecules $\sigma_i$ and $J$-molecules $\tau_i$, $i \in m$, such that

$$x_i = \varphi_i * \sigma_i, \quad y_i = \varphi_i * \tau_i, \quad C(\bar{\sigma}) = C(\bar{\tau}) \tag{*}$$

where $x_i, y_i$ are elements covered by pebbles $i$ in $\bar{I}, \bar{J}$ respectively. (Without loss of generality, we assume that, in the beginning, all $2m$ pebbles are on the board with all $x_i = y_i = \emptyset$.)

First we check that Duplicator can always play in the required manner. Clearly (*) holds in the initial position, as we can take all $\varphi_i = \emptyset$, all $\sigma_i = \sigma_j$ and all $\tau_i = \tau_j$. Now suppose that, after some number of steps, (*) holds witnessed by $\bar{\varphi}, \bar{\sigma}, \bar{\tau}$, and then Spoiler moves. By symmetry, we may assume that Spoiler moves pebble 0 on $\bar{I}$ from $x_0$ to $x_0'$. By Lemma 16, $x_0' = \varphi_0' *_I \sigma_0$ for some form $\varphi_0'$ and some molecule $\sigma_0'$ over $\bar{I}$. By Lemma 14, there exists a molecule $\tau_0'$ over $\bar{J}$ such that $C(\sigma_0', \sigma_1, \ldots, \sigma_{m-1}) = C(\tau_0', \tau_1, \ldots, \tau_{m-1})$. Duplicator can move pebble 0 on $\bar{J}$ from $y_0$ to $y_0' = \varphi_0' *_J \tau_0'$ and (*) will be restored.

Second, we check that (*) ensures that the map $x_i \mapsto y_i$ is a partial isomorphism and thus the proposed strategy of Duplicator is winning. For each $i, j \in m$, (*) implies $C(\sigma_i, \sigma_j) = C(\tau_i, \tau_j)$. By Lemma 17,

$$\begin{aligned}
x_i \in x_j &\iff \varphi_i * \sigma_i \in \varphi_j * \sigma_j \iff \text{In}(\varphi_i, \varphi_j, C(\sigma_i, \sigma_j)) \\
&\iff \text{In}(\varphi_i, \varphi_j, C(\tau_i, \tau_j)) \iff \varphi_i * \tau_i \in \varphi_j * \tau_j \\
&\iff y_i \in y_j
\end{aligned}$$

and similarly with $=$ and Eq in place of $\in$ and In.

The Equivalence Theorem is proved. $\square$

## 9.5 A Generalization

Until now we have considered the case when the input vocabulary $\Upsilon_0$ is empty. Now we consider the case when $\Upsilon_0$ consists of unary predicates, say $P_0, \ldots, P_{c-1}$. Restrict attention to input structures where the $c$ basic relations are pairwise disjoint; recall that such input structures are called colored set with colors $P_0, \ldots, P_{c-1}$. An automorphism of a colored set $I$ is simply a color preserving permutation of the elements of $I$. Recall that $\bar{I}_k$ is the collection of $k$-symmetric elements of $\text{HF}(I)$ as well as the corresponding structure of vocabulary $\{P_0, \ldots, P_{c-1}, \in, \emptyset\}$. We shall indicate how to modify the proof of the Equivalence Theorem to obtain the following version of it for colored sets.

**Corollary 5** *Fix positive integers $c, k, m$. If $I$ and $J$ are colored sets, in each of which all the colors $P_0, \ldots, P_{c-1}$ are sufficiently large, then $\bar{I}_k$ and $\bar{J}_k$ are $L^m_{\infty,\omega}$-equivalent.*

**Proof** To prove this, we need only to make the following changes in the proof of the Equivalence Theorem for naked sets.

First, a configuration should specify not only how the ranges of molecules overlap but also the colors of the atoms in the molecules. Thus, the configuration $C(\bar{\sigma})$ of a finite sequence $\bar{\sigma} = (\sigma_0, \ldots, \sigma_{l-1})$ of molecules should be defined as a pair $(C_=(\bar{\sigma}), C_*(\bar{\sigma}))$ where $C_=(\bar{\sigma})$ is the equivalence relation on $l \times k$ that we previously called the configuration and where $C_*(\bar{\sigma})$ is the function $l \times k \longrightarrow c$ sending each pair $(i, p)$ to the unique $r$ with $\sigma_i(p) \in P_r$. An abstract $l$-ary configuration is a pair $E$ whose first component $E_=$ is what we previously called an abstract $l$-ary configuration and whose second component $E_*$ is a function from $l \times k$ into $c$ that is constant on every equivalence class of $E_=$.

Next, we check that Lemma 14 still holds with this new notion of configuration. Two things must be added to the earlier proof of the lemma: $\tau'_0$ and $\sigma_0$ agree as to colors, and $\tau'_0$ can be extended to $\tau_0$ so as to maintain agreement with $\sigma_0$. The latter is clear because the colors $P_r$ in our input structures are large enough. As for the former, we must prove that, if $i \geq 1$ and $(0, p)Q_=(i, q)$ then $\tau_i(q)$ has the same color as $\sigma_0(p)$. But from $(0, p)Q_=(i, q)$ we get $\sigma_0(p) = \sigma_i(q)$, and this element has the same color as $\tau_i(q)$ because $C(\sigma_1, \ldots, \sigma_l)_* = C(\tau_1, \ldots, \tau_l)_*$.

In the statement of Lemma 15, "permutation" must be changed to "automorphism". The proof of that lemma is unchanged except that the computation verifying that $C(\pi^{-1}\rho, \sigma) = C(\rho, \pi\sigma)$ now verifies only that the $C_=$ components of the configurations agree. To get agreement of the $C_*$ components, we use the fact that $\pi$ is an automorphism and thus preserves colors.

The only other change in the earlier proof occurs in Lemma 16, where the difficult direction involved constructing a certain permutation $\pi$. In the colored situation, we must make sure that $\pi$ is an automorphism. For this purpose, we first verify that the $\pi_0$ defined in the earlier proof preserves colors. For atoms $a$ with $\pi_0(a) = a$, this is trivial, so we consider an atom $a$ for which $a = \tau_y(p)$ and $\pi_0(a) = \rho(p)$. Because $C(\tau_y, \sigma) = C(\rho, si)$, in particular the $C_*$ components agree. So $\tau_y(a)$ has the same color as $\rho(p)$. But these atoms are $a$ and $\pi_0(a)$, so $\pi_0$ preserves the color of $a$. Finally, when extending the map $\pi_0$ to an automorphism $\pi$, we must choose the extension so as to preserve colors, but this is trivially possible. □

# 10 Negative Results

## 10.1 Parity

Recall that a naked set is an input structure of the empty vocabulary.

**Theorem 10 (Parity Theorem)** *Parity is not in $\tilde{C}PTime$. Moreover, suppose that $K_1, K_2$ are disjoint infinite classes of naked structures, each containing structures of infinitely many cardinalities; then $(K_1, K_2)$ is not $\tilde{C}PTime$.*

**Proof** Let $\bar{\Pi}$ be any PTime ASM program with empty input vocabulary. We show that there are $I_1 \in K_1$ and $I_2 \in K_2$ such that $\bar{\Pi}$ does not distinguish between $I_1$ and $I_2$. By the Support Theorem in Section 8, there is a positive integer $k$ such that, in every run of $\bar{\Pi}$, every active set has a support of size $\leq k$.

Let Active($I$) be as in Section 6. By the First Fixed-Point Theorem in Section 6, there exists an FO+LFP sentence $\varphi$ that asserts that $\bar{\Pi}$ accepts $I$. (The sentence $\varphi$ asserts that there exist $i$ such that $D_f(i, true)$ holds in Active($I$) where $f$ is Output.) By Proposition 2 in Section 2, there is $m$ such that $\varphi$ is expressible in $L^m_{\infty,\omega}$. By the Equivalence Theorem in Section 9, $\varphi$ does not distinguish between any sufficiently large input structures $I_1, I_2$. □

## 10.2 Bipartite Matching is not in Choiceless PTime

Bipartite Matching is the following decision problem.

**Instance:** A bipartite undirected graph $G = (V, E)$ with the two parts (of boys and girls respectively) of the same size.

**Question:** Does there exists a perfect matching for $G$?

Recall that a *perfect matching* is a set $F$ of edges such that every vertex is incident to exactly one edge in $F$. A *partial matching* is an edge set $F$ such that every vertex is incident to at most one edge in $F$. The standard perfect-matching algorithm starts with the empty partial matching $F$ and then enlarges $F$ in a number of iterations. During each iteration, one constructs an auxiliary set $D$ of directed edges, then seeks a $D$-path $P$ from an unmatched boy to an unmatched girl, and then modifies $F$ by means of $P$.

We use variables $b, g$ to vary over boys and girls respectively. If $X$ is a set of edges, let

$$\begin{aligned}
\text{Boys-to-girls}(X) &= \{(b,g) : \{b,g\} \in X\} \\
\text{Girls-to-boys}(X) &= \{(g,b) : \{b,g\} \in X\}
\end{aligned}$$

And if $X$ is a set of ordered pairs of the form $(b,g)$ or $(g,b)$, let

$$\text{Unordered}(X) = \{\{b,g\} \,:\, (b,g) \in X \vee (g,b) \in X\}.$$

In Table 10.2, we give a self-explanatory program in the ASM language with the the choice construct for the perfect matching algorithm. It is customary to omit the keywords `do-in-parallel/enddo`. For readability, we take some little additional liberties with the ASM syntax.

The relation REACHABLE and the function PATH in the last transition rule are *external*. In other words, we take for granted algorithms that, given a boy $b$, a girl $g$ and set $D$ of directed edges over $V$, check whether there exists a $D$-path from $b$ to $g$ and if yes then construct such a path. For simplicity, we identify a path with the set of its edges.

For the benefit of those unfamiliar with the algorithm, let us explain one iteration of the algorithm in the case when the given bipartite graph has a perfect matching $M$. Suppose that $F$ is a current partial matching involving $k$ boys, and there are some $F$-unmatched boys. Abusing notation, let $M(g)$ be the girl $M$-matched to $b$ and let $F(g)$ be the boy $F$-matched to $g$. Let $D$ be as in the Build-Digraph rule. For any $F$-ummatched boy $b_0$, there exists a $D$-path $P$ from $b_0$ to some $F$-unmatched girl. Indeed, construct

$$\begin{aligned} g_1 &= M(b_0), b_1 = F(g_1) \\ g_2 &= M(b_1), b_2 = F(g_2) \\ &\ldots \\ g_k &= M(b_{k-1}), b_k = F(g_k) \\ g_{k+1} &= M(b_k) \end{aligned}$$

until you encounter an $F$-unmatched girl $g_{k+1}$. It is easy to see that, if $X = \text{Unordered}(P)$, then $(F - X) \cup (X - F)$ is a partial matching involving $k+1$ boys.

**Theorem 11 (Bipartite Matching Theorem)**
*Bipartite Matching is not in $\tilde{C}PTime$.*

**Proof** Given an even integer $n = 2p > 2$, we construct two bipartite graphs $G_0$ and $G_1$ on a set

$$V_n = \{b_0, \ldots, b_{n-1}\} \cup \{g_0, \ldots, g_{n-1}\}$$

of $n$ boys and $n$ girls. In $G_0$, (1) the first $p$ boys and the first $p$ girls form a complete bipartite graph, (2) the last $p$ boys and the last $p$ girls form a complete bipartite

```
if Mode = Initial then
  F := ∅, Mode := Examine
endif

if Mode = Examine then
  if there is an unmatched boy then
    Mode := Build-Digraph
  else Output := true, Mode := Final
endif

if Mode = Build-Digraph then
  D := Girls-to-boys(F) ∪ Boys-to-girls(E-F)
  Mode := Build-Path
endif

if Mode = Build-Path then
  choose an unmatched boy b
    if (∃ unmatched girl g) REACHABLE_D(b,g) then
      choose an unmatched girl g with REACHABLE(D,b,g)
        P := PATH_D(b,g), Mode := Modify-Matching
      endchoose
    else Output := false, Mode := Final
  endchoose
endif

if Mode = Modify-Matching then
  F := (F - Unordered(P)) ∪ (Unordered(P) - F)
  Mode := Examine
endif
```

Table 1: The perfect matching algorithm

graph, and (3) there are no other edges. Clearly, $G_0$ has a perfect matching. In $G_1$, (1) the first $p+1$ boys and the first $p$ girls form a complete bipartite graph, (2) the last $p-1$ boys and the last $p$ girls form a complete bipartite graph, and (3) there are no other edges. Clearly, $G_1$ has no perfect matching.

Notice that the two graphs are essentially 4-colored sets; "adjacency" is definable from the colors. The rest of the proof is similar to the proof of the Parity Theorem, except that, instead of the Support Theorem and Equivalence Theorem, we use Corollary 3 and Corollary 5 respectively. □

# 11 Revised Choiceless PTime

Enrich the computation model with a static nullary function InputSize. In the definition of initial states with $n$ atoms require that InputSize is the von Neumann ordinal for $n$. That changes the notion of $\tilde{C}$PTime; indeed, there is an obvious PTime ASM program with InputSize that accepts all naked sets of odd cardinality and rejects all naked sets of even cardinality. Let us call the new complexity class $\tilde{C}$PTime$^+$.

$\tilde{C}$PTime$^+$ is still defined in terms of pairs $(K_1, K_2)$ of classes of finite structures of the same vocabulary which are disjoint but not necessarily complementary. But the presence of InputSice allows us to augment each program with a polynomial clock and thus ensure that every input is either accepted or rejected within polynomial time. Thus we can return to the usual framework where a complexity class is defined in terms of classes of structures (rather than pairs of classes).

**Theorem 12** *There is a class of "standard" ASM programs with InputSize and the following properties:*

1. *Each standard program accepts or rejects every input within polynomial time.*

2. *Every PTime program can be simulated by a standard program.*

3. *There exists a PTime Turing machine that, given an ASM program, decides wether the program is or is not standard.*

**Proof** A standard program is an arbitrary program $\Pi$ augmented with a clock that bounds the number of microsteps in the computation of $\Pi$ by some polynomial $p$ of the InputSize. The clock is initialized to zero, is incremented at each microstep of $\Pi$, and stops the computation when it reaches $p$(InputSize). To satisfy the third claim, use a standard format for the clock program. □

Fix a particular class of standard programs and redefine $\tilde{C}$PTime$^+$ as follows: A class $K$ of finite structures of the same vocabulary is $\tilde{C}$PTime$^+$ if there exists a standard program that accepts exactly the structures in $K$.

**Corollary 6** *There is a PTime (two-valued) logic that expresses exactly $\tilde{C}PTime^+$ properties.*

**Proof** View standard programs as formulas. □

In the rest of this section, we show that:

- for input structures with just one unary relation $U$, the parity of $|U|$ cannot be computed in $\tilde{C}PTime^+$, and

- Bipartite Matching is outside $\tilde{C}PTime^+$.

For this purpose, we first observe that Corollary 3, the Support Theorem for colored structures, remains true when InputSize is allowed, as a static nullary function, in programs. Indeed since the value of InputSize is a von Neumann ordinal, it involves no atoms and is therefore fixed by all permutations of the atoms. It follows that, when any program $\Pi$ is run with a colored set $I$ as the input structure, the set $\mathcal{A}$ of active elements is invariant under all automorphisms of $I$. With this observation, the proof of the Support Theorem for colored structures goes through as before.

Recall from Section 9 that, for any colored set $I$ and any positive integer $k$, $\bar{I}_k$ denotes the set of $k$-symmetric elements of $\mathrm{HF}(I)$ as well as the corresponding structure with $\in$, $\emptyset$, and the colors. For our present purposes, we must consider the expansion $(\bar{I}_k, |I|)$ of the structure $\bar{I}_k$ where the cardinality $|I|$ of the input set is named by the nullary symbol InputSize. We show next that the only effect of this extra constant on the equivalence theorem is to restrict the result (as one might expect) to input structures of equal cardinality.

**Theorem 13** *Fix positive integers $c$, $k$, and $m$. If $I$ and $J$ are $c$-colored sets of the same cardinality, and if all the colors $P_0$, ..., $P_{c-1}$ are sufficiently large in both of them, then $(\bar{I}_k, |I|)$ and $(\bar{J}_k, |J|)$ are $L_{\infty,\omega}^m$-equivalent.*

**Proof** Observe first that, for each natural number $n$, there is a form $\varphi_n$ such that $\varphi_n * \sigma = n$ for every molecule $\sigma$. Such $\varphi_n$ can be defined inductively by

$$\varphi_n = \{(\varphi_r, E) \mid r < n \text{ and } E \text{ is a binary configuration}\}.$$

The verification that $\varphi_n * \sigma = n$ is a trivial induction on $n$.

Now suppose $I$ and $J$ are as in the hypothesis of the theorem. The $m$-pebble game for the structures $(\bar{I}_k, |I|)$ and $(\bar{J}_k, |J|)$ is the same as the $(m+1)$-pebble game for $\bar{I}_k$ and $\bar{J}_k$ with one pebble located permanently at the natural number $|I| = |J| = n$ in both structures. Since this number is the denotation of the same form $\varphi_n$ in both structures, Duplicator can still use the winning strategy described in Section 9: match the forms and the configurations of molecules. □

Let Subset Parity be the following decision problem.

*Instance:* A structure $(I, U)$ where $U$ is a unary relation on $I$ (i.e., $U \subseteq I$).

*Question:* Is $|U|$ odd?

**Corollary 7** *Subset Parity is not in $\tilde{C}PTime^+$.*

**Proof** If Subset Parity were in $\tilde{C}PTime^+$, then, by Proposition 2, Theorem 4, and Corollary 3, there would be positive integers $m$ and $k$ such that, whenever $(I, U)$ is a positive instance and $(J, V)$ a negative instance of Subset Parity, then Duplicator has no winning strategy in the $m$-pebble game for $(\bar{I}_k, |I|)$ and $(\bar{J}_k, |J|)$ (where $I$ abbreviates $(I, U)$ and similarly for $J$).

On the other hand, structures of the form $(I, U)$ can be regarded as 2-colored sets, with colors $U$ and $I - U$. Thus, by Theorem **??**, Duplicator has a winning strategy provided $|I| = |J|$ and all of $|U|, |V|, |I - U|$, and $|J - V|$ are large enough. This situation is clearly compatible with $|U|$ being odd and $|V|$ even, so we have a contradiction. □

**Corollary 8** *Bipartite Matching is not in $\tilde{C}PTime^+$.*

**Proof** We can use exactly the same proof as in Section 10, because the two structures used there had the same cardinality. □